\documentclass[twoside,12pt]{article}
\usepackage{latexsym,newlfont,amsfonts,amssymb,makeidx}

\newcounter{gavriil}
\usepackage{amsmath}

\DeclareMathOperator{\Alt}{Alt}

\DeclareMathOperator{\ad}{ad}
\DeclareMathOperator{\B}{B}
\DeclareMathOperator{\C}{C}

\DeclareMathOperator{\differential}{d}

\DeclareMathOperator{\End}{End}
\DeclareMathOperator{\Ha}{H}
\DeclareMathOperator{\Hom}{Hom}
\DeclareMathOperator{\Id}{Id}
\DeclareMathOperator{\id}{id}
\DeclareMathOperator{\im}{Im}

\DeclareMathOperator{\U}{U}
\DeclareMathOperator{\Vect}{Vect}
\DeclareMathOperator{\Z}{Z}

\newenvironment{proof}
{\medskip\noindent {\sc Proof:}} {\qed\medskip}

\newcommand{\nm}{\small\sc}

\newcommand{\cL}{\mathcal{L}}

\newcommand{\cP}{\mathcal{P}}

\newcommand{\ba}{\mathbf{a}}

\newcommand{\bb}{\mathbf{b}}

\newcommand{\bl}{\mathbf{l}}

\newcommand{\br}{\mathbf{r}}

\newcommand{\bu}{\mathbf{u}}

\newcommand{\bw}{\mathbf{w}}

\newcommand{\al}{\alpha}
\newcommand{\be}{\beta}

\newcommand{\di}{\differential\!}
\newcommand{\De}{\Delta}
\newcommand{\de}{\ensuremath{\delta}}

\newcommand{\ep}{\varepsilon}

\newcommand{\Ga}{\Gamma}
\newcommand{\g}{\mathfrak g}
\newcommand{\ga}{\gamma}
\newcommand{\gd}{{\mathfrak{g}}^*}
\newcommand{\geqs}{\geqslant}

\newcommand{\hb}{\ensuremath{\hbar}}

\newcommand{\ka}{\ensuremath{\kappa}}
\newcommand{\kk}{\mathbb{C}}

\newcommand{\La}{\ensuremath{\Lambda}}
\newcommand{\la}{\ensuremath{\lambda}}

\newcommand{\leqs}{\leqslant}

\newcommand{\nl}{\newline}
\newcommand{\ol}{\overline}
\newcommand{\Omg}{\Omega}
\newcommand{\omg}{\omega}

\newcommand{\ot}{\otimes}
\newcommand{\Ph}{\Phi}
\newcommand{\ph}{\varphi}
\newcommand{\qed}{\ensuremath{\blacksquare}}
\newcommand{\rk}{\mathbb R}

\newcommand{\si}{\sigma}
\newcommand{\T}{\mathrm T}

\newcommand{\tht}{\theta}
\renewcommand{\to}{\longrightarrow}
\newcommand{\Ug}{{\U}(\g)}
\newcommand{\Ugh}{{\U}_\hb(\g)}

\newcommand{\un}{{\mathbf 1}}  
\newcommand{\ups}{\upsilon}
\newcommand{\w}{\wedge}
\newcommand{\wtl}{\widetilde}
\newcommand{\z}{\mathbb Z}
\newcommand{\ze}{\zeta}

\makeindex

\newtheorem{thm}{Theorem}[section]
\newtheorem{prop}{Proposition}[section]
\newtheorem{lm}{Lemma}[section]
\newtheorem{cor}{Corollary}[section]

\newtheorem{defn}{Definition}[section]
\newtheorem{eg}{Example}[section]
\newtheorem{rem}{Remark}[section]

\newcommand{\Ah}{A_\hb}
\renewcommand{\a}{\mathfrak{a}}
\renewcommand{\bw}{\bigwedge}
\newcommand{\Caut}{{\tht}}
\newcommand{\Car}{\mathfrak{h}}
\newcommand{\clal}{{\ol{\al}}}
\newcommand{\clbe}[1][]{{\ol{\be}_{#1}}}
\newcommand{\comp}{\circ}
\newcommand{\Db}{\De_\hb}
\newcommand{\dhb}{[\![\hb]\!]}
\newcommand{\dl}{[\![}
\newcommand{\dr}{]\!]}
\newcommand{\dSum}{\bigoplus} 
\newcommand{\ds}{\wtl{\di}}
\newcommand{\dsum}{\oplus} 
\newcommand{\Fb}{F_\hb}
\newcommand{\fAh}{A\dhb}

\newcommand{\gut}{admissible}
\renewcommand{\k}{\mathfrak{k}}
\newcommand{\Mal}[1][l]{M_{#1\al}}
\newcommand{\m}{\mathfrak{m}}
\newcommand{\mal}[1][]{\m_{#1\ol{\al}}}

\newcommand{\mdh}[1][n]
    {\ \ \left(\hspace{-12pt}\mod\hb^{#1}\right)}
\newcommand{\mdhs}[1][n]
    {\ \ \left(\hspace{-8pt}\mod\hb^{#1}\right)}
\newcommand{\mh}{\mu_\hb}
\newcommand{\mser}[1][]{\mu_0 + \mu^{#1}_1\hb 
            + \mu^{#1}_2\hb^2 + \mu^{#1}_3\hb^3 + \cdots}

\newcommand{\Pb}{\Ph_\hb}
\newcommand{\rf}{\nm}
\newcommand{\tmh}{\wtl{\mu}_\hb}
\newcommand{\Ughr}{\U_\hb(\g,r)}
\newcommand{\vf}{\mathfrak{v}}
\newcommand{\vs}{\vspace{.2in}}
\newcommand{\Wal}[1][l]{\Wurz_{#1\al}}

\newcommand{\UWurz}{\mathrm{P}}
\newcommand{\Wurz}{\Omg}

\begin{document}

\bibliographystyle{plain}

\begin{titlepage}
\thispagestyle{empty}
\title{$\Ughr$ Invariant Quantization\\ 
on Some Homogeneous Manifolds}
\author{Vadim Ostapenko}
\date{\ }
\maketitle

\thispagestyle{empty}
\vspace{1in}
\begin{abstract}
We consider a class of homogeneous manifolds over a simple {\nm Lie}
group which appears in the problem of classification of homogeneous 
manifolds with reductive subgroups of maximal rank as stabilizer of
a point.
We prove that any manifold of this class possesses a
Poisson bracket admitting a quantization invariant
(equivariant) with
respect to the corresponding quantum group.
\end{abstract}

\newpage
\thispagestyle{empty}
     \ 
\end{titlepage}

\thispagestyle{empty}
\tableofcontents

\ \ 
\newpage
\thispagestyle{empty}
\ \ \ 
\newpage

\section*{Acknowledgments}

I am deeply grateful to my scientific advisers {\nm Joseph F. Donin}
and {\nm Steven Shnider} for posing the problem and effective 
guidance.
I thank to my friends {\nm Vladimir G. Berkovich} and 
{\nm  Vladimir A. Hinich} who taught me homological algebra and many 
other things in mathematics.
I thank also {\nm Arkadi\u{i} L. Onishchik} and 
{\nm Ernst B. Vinberg} for very helpful discussions concerning roots 
systems.
In particular, {\nm E.~B.~Vinberg} pointed out to the second 
claim of Lemma~\ref{th:mBarAlIrredInTensProd}.
I thank {\nm Alexander A. Stolin} for useful discussing the
Belavin---Drinfeld r-matrices.

I thank my wife {\nm Jenny} for support she gave me during doing the 
research.
My special thanks to my parents {\nm Ninel} and {\nm Borislav 
Ostapenko} and to my parents-in-law {\nm Marita} and 
{\nm Isaac Gun} who, spending a lot of time with our children, 
made it possible for me to carry out this research.

Finally, I wish to express my gratitude to the Department of 
Mathematics and Computer Science of Bar-Ilan University, 
Ramat Gan, for the enthusiastic atmosphere created both by  
faculty members and numerous outstanding visiting mathematicians.

\newpage
\ \ 
\newpage

\section{Introduction}

A quantum homogeneous manifold is obtained from the usual
homogeneous manifold by replacing the original commutative 
function algebra with a deformed non-commutative algebra.
Since a homogeneous manifold is equipped with a {\nm Lie} group 
action, it is natural to look for deformation quantizations of the 
function algebra which are invariant with respect to the action of 
the corresponding quantum group.

Let $G$ be a simple connected {\nm Lie} group over $\kk$, $\g$ its 
{\nm Lie} algebra, $K$ a closed subgroup of $G$, $\k$ the {\nm Lie}
algebra of $K$.
Denote by $M$ the homogeneous manifold $G/K$, by $\mu$ the 
commutative multiplication in $\C^\infty(M)$.
A quantization of $M$ is a formal deformation
$\mh=\mu+\hb\mu_1+\hb^2\mu_2+\cdots$ of $\mu$ which defines an
associative multiplication in the space $\C^\infty\dhb$ of formal 
power series in $\hb$ with smooth functions on $M$ as their 
coefficients.
Let $\Ugh=(\Ug,\Db)$ be a quantization of $\g$.
We say that $\mh$ is invariant with respect to $\Ugh$ if 
$x.\mh(a,b)=\mh\left(\Db(x)\cdot(a\ot b)\right)$ for any 
$a,b\in\C^\infty(M)$ and $x\in\Ugh$.
\index{multiplication!invariant}
One can assume from the beginning that the 
bilinear mapping $\mu_1$ is skew symmetric.
The associativity of the deformed multiplication $\mh$ implies
the {\nm Jacobi} identity for $\mu_1$.
Thus $\mu_1$ is essentially a {\nm Poisson} bracket on $M$.
The first question we investigate is whether there exist 
{\nm Poisson} brackets on the homogeneous manifold $M=G/K$
admissible for $\Ugh$ invariant quantization.
The second question we consider is whether there exists a deformed 
multiplication $\mh$ which is invariant under the quantum group
action.

It turns out that brackets which admit invariant quantization are 
of special form which we describe now.
The $G$ action on $M$ determines a homomorphism from $\g$ to the 
{\nm Lie} algebra $\Vect(M)$\index{$\Vect(M)$} of vector fields on 
$M$ and so it induces a linear map $\rho:\bw^2\g\to\bw^2\Vect(M)$. 
Consider a bivector $r\in\bw^2\g$ such that the {\nm Schouten} 
bracket $\dl\rho(r),\rho(r)\dr$ is $\g$ invariant.
Such an element is called a {\nm Belavin--Drinfeld} classical 
r-matrix.
Since the algebra $\g$ is simple, there exists a $\g$ invariant 
$3$-vector $\ph\in\bw^3\g$ unique up to constant multiple.
We normalize $r$ so that $\dl r,r\dr=\ph$ and thus 
$\dl\rho(r),\rho(r)\dr=\rho(\ph)$.
Each {\nm Belavin--Drinfeld} classical r-matrix $r$ defines a 
{\nm Lie} bialgebra structure on $\g$.
{\rf P.~Etingof} and {\rf D.~Kazhdan} have proven
\cite{EtinhoffKazdan:QuantBialg} that any {\nm Lie} bialgebra can
be quantized.
In particular, if such a structure is determined by $r$, there
exists a quantum group $\Ughr$ whose multiplication is preserved 
from $\Ug$ and the comultiplication is of the form 
$\Db=\De+\hb r+\cdots$ where $\De$ denotes the original 
comultiplication in $\Ug$.
If $\rho(\ph)=0$ then $\rho(r)$ is a {\nm Poisson} bracket on $M$
which is called an r-matrix {\nm Poisson} bracket.
It can happen, however, that $\rho(\ph)\not=0$ on $M$, but there
exists a $G$ invariant bivector field $s$ on $M$ such that
\begin{equation}\label{eq:mCYBE}
\dl s,s\dr=-\rho(\ph)
\end{equation}
We call such $s$ a 
\emph{$\ph$-Poisson bracket}\index{P{\nm oisson} bracket} on $M$.
Thus $\ph$-{\nm Poisson} bracket is a skew symmetric bracket obeying 
the {\nm Leibniz} rule and the weak version of the {\nm Jacobi} 
identity expressed by equation (\ref{eq:mCYBE}).

For a $\ph$-{\nm Poisson} bracket $s$, the sum $s+\rho(r)$ is a 
{\nm Poisson} bracket on $M$, but it is not $\g$ invariant.
We call {\nm Poisson} brackets of this form 
\emph{\gut}\index{P{\nm oisson} bracket!\gut}.
{\rf J.~Donin}, {\rf D.~Gurevich} and {\rf S.~Shnider} proved
\cite{DoninGurevicShnider:DoubleQuantiz} that if {\nm Poisson} 
bracket on a homogeneous manifold admits a $\Ughr$ invariant
quantization, it is of necessity \gut,
$s$ and $r$ satisfying (\ref{eq:mCYBE}).
If such a {\nm Poisson} bracket exists on $M$, the question 
arises whether it can be quantized in such a way that the deformed 
multiplication $\mh=\mu+\hb s+\hb^2\mu_2+\cdots$ is $\Ughr$ 
invariant.
We call an \gut\ {\nm Poisson} bracket 
\emph{quantizable}\index{P{\nm oisson} bracket!quantizable} if
there exists some its $\Ughr$ invariant quantization.

The problem of invariant quantization of {\rf Poisson} brackets was 
considered by many authors, among them are {\rf J.~Donin}, 
{\rf D.~Gurevich}, {\nm S.~Khoroshkin},
{\rf Sh.~Majid},  {\rf A.~Radul}, {\rf V.~Rubtsov}, {\rf S.~Shnider}
and others.

{\rf J.~Donin, D.~Gurevich} and {\rf Sh.~Majid} have considered 
\cite{DoninGurevicMajid:RMatrBrackets} so called 
{\nm Drinfeld--Jimbo} classical r-matrix which is a particular case 
of {\nm Belavin--Drinfeld} r-matrix.
They proved that if the {\nm Lie} subalgebra $\k$ contains a
maximal nilpotent subalgebra, then the {\nm Drinfeld--Jimbo}
classical r-matrix generates a {\nm Poisson} bracket on $M=G/K$,
and that there exists a quantization of this bracket invariant with
respect to the {\nm Drinfeld--Jimbo} quantum group $\Ugh$.
{\rf J.~Donin} and {\rf D.~Gurevich} 
\cite{DoninGurevic:DrJimbPoissBr} proved that on a semi-simple orbit 
of coadjoint representation the {\nm Poisson} bracket coming from 
the {\nm Sklyanin--Drinfeld} bracket on $G$ is quantizable.
{\rf J.~Donin} and {\rf S.~Shnider} \cite{DoninShnider:QSymSpaces} 
proved that the {\nm Drinfeld--Jimbo} classical r-matrix 
generates a {\nm Poisson} bracket on any symmetric space, and 
solved the problem of quantization for this bracket.

As we mentioned above, {\rf J.~Donin}, {\rf D.~Gurevich} and 
{\rf S.~Shnider} proved \cite{DoninGurevicShnider:DoubleQuantiz} 
that if {\nm Poisson} bracket on a homogeneous manifold admits a 
$\Ughr$ invariant quantization, it is of necessity of the form 
$s+\rho(r)$ with $s$ and $r$ satisfying (\ref{eq:mCYBE}).
They have shown that almost any such bracket can be quantized 
invariantly on any semi-simple orbit of the coadjoint 
representation.
In the present work, we introduce another class of homogeneous 
manifolds containing manifolds which are not necessarily orbits of 
the coadjoint representation of $\g$.
We prove that any manifold in of this class can be equipped with an
essentially unique {\nm Poisson} bracket which admits a $\Ughr$ 
invariant quantization.

To describe this class of homogeneous manifolds, fix a {\nm Cartan}
subalgebra $\Car$ in $\g$, a simple root base for the corresponding
root system and an integer $l\geqs 2$.
Take a simple root $\al$ and denote by $\k$ the {\nm Lie} 
subalgebra of $\g$ generated by the {\nm Cartan} subalgebra and by 
all roots for which the coefficient of $\al$ is divisible by $l$.
Denote by $K$ the subgroup of $G$ corresponding to $\k$ and set 
$\Mal=G/K$.
The significance of the manifolds $\Mal$ is in the fact that any
quotient of $G$ by a reductive subgroup of maximal rank can be 
obtained by taking consequent quotients of homogeneous manifolds 
of this type, see \cite{DynkinE:SSSubAlgsInSSLieAlgs}.

We prove that for any $\Mal$ there exists a $\g$ invariant bivector 
field $s$ satisfying (\ref{eq:mCYBE}).
We prove then that for any {\nm Belavin--Drinfeld} r-matrix $r$, 
the {\nm Poisson} bracket $\rho(r)+s$ can be quantized in 
$\Ughr$ invariant way.
We produce this quantization in two steps.
First, using methods developed from the techniques of 
\cite{DoninShnider:QSymSpaces}, we prove that the $\ph$-{\nm Poisson}
bracket $s$ can be quantized in such a way that the deformed 
multiplication is invariant under the action of the group $G$ and 
obeys some deformed associativity constraint.
Then, using $\rho(r)$, we correct the above deformed multiplication, 
turning it into a new multiplication satisfying the usual 
associativity law.
This new multiplication is invariant under the action of the quantum
group $\Ughr$.
The reason why we first pass to the category with non-trivial
associativity is that while losing the associativity, we gain $G$
invariance. 
The next step, passing to the quantum group symmetry, is 
achieved by an equivalence of categories.
To construct this equivalence, we prove that there exists an
invertible element $\Fb\in\Ug^{\ot 2}\dhb$ such that\index{$\Fb$}
\begin{equation}\label{eq:ThEtinhKazhd}
\Db(x)=\Fb^{-1}\cdot\De(x)\cdot\Fb
\end{equation}
for any $x\in\Ug$.
Here $\De$ is the standard comultiplication in $\Ug$ and $\Db$ is 
the comultiplication in the quantum group $\Ughr$.
This proof is essentially based on the work 
\cite{EtinhoffKazdan:QuantBialg} by {\rf P.~Etingof} and 
{\rf D.~Kazhdan}.
This two step method allows us to reduce $\Ughr$ invariant
quantizations of brackets $s+\rho(r)$ for different 
{\nm Belavin--Drinfeld} r-matrix $r$ to a single quantization of 
the $\ph$-{\nm Poisson} bracket $s$ in the category with non-trivial
associativity.

The work is organized as follows.
In Section~\ref{sect:InfinitesimalsForAlgDeforms} we recall some 
facts about {\nm Hochschild} complexes, with special attention to 
those aspects which are important for our purposes.
In particular, we consider some symmetry properties of  
{\nm Hochschild} complexes.

In Section~\ref{sect:MonCategDrAlgs} we explain how to reduce the
problem of associative $\Ughr$ invariant quantization to the 
problem of non-associative $G$ invariant quantization.
Namely, for any {\nm Belavin---Drinfeld} r-matrix, we construct the
element $\Fb$ which follows (\ref{eq:ThEtinhKazhd})
and gives the equivalence of the categories with trivial and 
non-trivial associativity.

In Section~\ref{sect:InfinitesForInvarQn} we study properties of 
infinitesimals for $G$ invariant quantization on homogeneous $G$ 
manifolds in the category with the non-trivial associativity.
We introduce the notion of $\ph$-{\nm Poisson} bracket and show 
that, similarly to the case of associative deformation, any $G$ 
invariant deformation with non-trivial associativity has a 
$\ph$-{\nm Poisson} bracket as its infinitesimal.

In Section~\ref{sect:PhiPoissQuantization} we prove a general
theorem which gives a sufficient condition for quantizability of
$\ph$-{\nm Poisson} brackets.

In Section~\ref{sect:Examples} we introduce the homogeneous 
manifolds $\Mal$.
We prove that all manifolds $\Mal$ posses $G$ invariant 
$\ph$-{\nm Poisson} brackets and give explicit form for all these 
brackets.
We consider cochain complexes generated by these brackets and 
prove that dimensions of some cohomologies is equal to zero.
This allows us to apply the general theorem of 
Section~\ref{sect:PhiPoissQuantization} and prove that these 
$\ph$-{\nm Poisson} brackets  can be quantized $G$ invariantly in the
 category with non-trivial associativity.
Applying the results of Section~\ref{sect:MonCategDrAlgs}, we
conclude that any {\nm Poisson} bracket $s+\rho(r)$ can be 
quantized invariantly with respect to the quantum group $\Ughr$ 
action.

\newpage
\section{Poisson brackets as infinitesimals\\ 
for algebra deformations}
\label{sect:InfinitesimalsForAlgDeforms}

In this section we give basic definitions of deformation theory 
of commutative algebras, with special attention to the function 
algebras on manifolds.
We show that the linear term of a formal deformation of commutative 
algebra is a {\nm Poisson} bracket on this algebra.
We recall some facts about {\nm Hochschild} complexes, stressing 
their symmetry properties.

Fix an associative commutative $\kk$ algebra $A$ with unity and
denote by $\mu$ the multiplication in $A$.
We denote by $\kk\dhb$ the $\kk$ algebra of formal power series in a 
variable $\hb$.
All tensor products over $\kk\dhb$ are assumed to be \emph{completed}
in the $\hb$-adic topology.

\subsection{Deformations of commutative algebras}

\begin{defn}\label{defn:AlgDfrms}
A $\kk\dhb$ algebra $(\Ah,\mh)$ is called 
\emph{a (formal) deformation}\index{deformation of algebra}
\index{algebra deformation} of $(A,\mu)$ if \nl
\mbox{\hspace{.15in}}\emph{(i)}\mbox{\hspace{.15in}}
 $\Ah$ is equal to $A\ot\kk\dhb$ as a topological $\kk\dhb$ module;\nl
\mbox{\hspace{.13in}}\emph{(ii)}\mbox{\hspace{.13in}} 
$\Ah/\hb \Ah\cong A$ as $\kk$ algebras.
\end{defn}\index{$\mh$}
Thus a deformation of $A$ is a topologically free $\kk\dhb$ 
module $\fAh$ equipped  with a $\kk\dhb$ linear mapping
$\mh: \fAh\ot_{\kk\dhb}\fAh\to\fAh$ obeying the associativity law.
One can think about $\mh$ as of a formal series $\mser$ where $\mu_i$ 
is the $\kk\dhb$ linear extension of a $\kk$ linear mapping 
$A\ot A\to A$.
Usually,  $\mu_0$ is called 
\emph{the initial term of $\mh$}\index{initial term of deformation}
and $\mu_1$ is called 
\emph{the infinitesimal of $\mh$}\index{infinitesimal of deformation}.
The second condition in the definition means that $\mu_0=\mu$, the 
original multiplication in $A$.
We are interested in a special kind of deformations which are called 
\emph{quantizations}\index{quantization}.
Before giving the definition, we examine some important properties 
of infinitesimals.

Let $(\Ah,\mh)$ be a deformation of $(A,\mu)$.
The associativity of $\mh=\mu+\hb\mu_1+\hb^2\mu_2+\cdots$ is
equivalent to the following infinite system of equations:
$$
\sum_{\begin{array}{c}
\scriptstyle
i+j=n\\
\scriptstyle 
 i,j\geqs 0
\end{array}}\hspace{-.1in}\Big(
\mu_i\comp(\mu_j\ot\id)-
\mu_i\comp(\id\ot\mu_j)\Big)
=0,\ \ \ \ n=0,1,2,\dots.
$$
Separating the terms with $i=0$ and $j=0$, one obtains the above 
equation in the following form:
\begin{eqnarray*}
\mu\comp(\id\ot\mu_n)-\mu_n\comp(\mu\ot\id)&+&
\mu_n\comp(\id\ot\mu)-\mu\comp(\mu_n\ot\id)=\\
&=&\hspace{-.1in}\sum_{\begin{array}{c}
\scriptstyle
i+j=n\\
\scriptstyle 
 i,j\geqs 1
\end{array}}
\hspace{-.1in}\Big(\mu_i\comp(\mu_j\ot\id)-
\mu_i\comp(\id\ot\mu_j)\Big).
\end{eqnarray*}
The expression in the left hand side is known as the 
{\nm Hoch\-schild} coboundary of $\mu_n$, and is denoted by 
$\di\mu_n$.
This gives an equivalent form for the associativity of $(\Ah,\mh)$:
\begin{equation}\label{eq:CartMaur}
\di\mu_n
=\hspace{-.1in}\sum_{\begin{array}{c}
\scriptstyle
i+j=n\\
\scriptstyle 
 i,j\geqs 1
\end{array}}
\hspace{-.1in}\Big(\mu_i\comp(\mu_j\ot\id)-
\mu_i\comp(\id\ot\mu_j)\Big).
\end{equation}
For $n=1$ one has $\di\mu_1=0$, i.e. the infinitesimal of any 
deformation is always a {\nm Hochschild} cocycle.
\begin{defn}\label{def:EquivDeforms}
Two deformations, $(A\dl\hb\dr,\mh)$ and $(A\dl\hb\dr,\mh')$, 
are {\em equivalent} if there exists  a $\kk\dl\hb\dr$ module 
automorphism $u_\hb$ of  $\fAh$ such that\nl
\mbox{\hspace{.15in}}\emph{(i)}\mbox{\hspace{.15in}} 
the restriction $u_0$ of $u_\hb$ to $A$ is the identity map;\nl
\mbox{\hspace{.13in}}\emph{(ii)}\mbox{\hspace{.13in}} 
the following diagram is commutative:

\begin{picture}(200,120)(-150,-15)
\put(15,70){\makebox(0,0){$\fAh\ot\fAh$}}
\put(115,70){\makebox(0,0){$\fAh$}}
\put(15,10){\makebox(0,0){$\fAh\ot\fAh$}}
\put(115,10){\makebox(0,0){$\fAh$}}
\put(53,70){\vector(1,0){40}}
\put(51,10){\vector(1,0){45}}
\put(70,80){\makebox(0,0){$\mh$}}
\put(70,20){\makebox(0,0){$\mh'$}}
\put(15,55){\vector(0,-1){30}}
\put(115,55){\vector(0,-1){30}}
\put(-10,40){\makebox(0,0){$u_\hb\ot u_\hb$}}
\put(127,40){\makebox(0,0){$u_\hb$}}
\end{picture}

\end{defn}

\subsection{The Hochschild complex}\label{subsect:HochCompl}

The obstruction theory, which links formal deformations of algebras
to {\nm Hoch\-schild} complexes, was developed by 
{\rf Murray Gerstenhaber} in 
\cite{GerstenhaberM:DeformsRingsAlgs1},
\cite{GerstenhaberM:DeformsRingsAlgs2},
\cite{GerstenhaberM:DeformsRingsAlgs3} and
\cite{GerstenhaberM:DeformsRingsAlgs4}.
We give here some properties of {\nm Hochschild} cocycles on
commutative algebras.
For a given commutative $\kk$ algebra $(A,\mu)$,
\emph{the Hochschild complex}\index{H{\nm ochschild}!complex} on $A$ 
with coefficients in $A$ is the graded $\kk$ vector space
$$
\C^\bullet(A;A)=\bigoplus_{p\geqs 0}\Hom(A^{\ot p},A)
$$
together with the coboundary operator $\di:\C^p(A;A)\to\C^{p+1}(A;A)$ 
defined by
\begin{equation}\label{eq:DefHochschCobound}
\di\xi=
\mu\comp(\id\ot\xi)+\sum_{k=1}^p(-1)^k
\xi\comp(\id^{\ot k-1}\ot\mu\ot\id^{\ot p-k})
+(-1)^{p+1}\mu\comp(\xi\ot\id).\ \ \ \ 
\end{equation}\index{H{\nm ochschild}!coboundary operator}

There are two useful operators on $\C^\bullet(A;A)$.
The first one we denote by $\tau$ and define as
$$
(\tau.\xi)(a_1,a_2,\dots,a_p)=
(-1)^\frac{p(p+1)}{2}\xi(a_p,a_{p-1},\dots,a_1)
$$
for $\xi\in\C^p(A;A)$.
\begin{prop}\label{th:TauIsComplexMorphism}
The map $\tau$ commutes with the coboundary operator on 
$\C^\bullet(A;A)$ and thus induces a cochain complex morphism,
which we denote also by $\tau$.
This morphism splits the complex into the direct sum of two 
sub-complexes,
$\C^\bullet(A;A)=\C^\bullet_+(A;A)\dsum\C^\bullet_-(A;A)$, where
$\C^\bullet_+(A;A)$ consists of all elements $\xi$ with
$\tau.\xi=\xi$ and $\C^\bullet_-(A;A)$ consists of all elements $\xi$
such that $\tau.\xi=-\xi$.
\end{prop}\index{$\C^\bullet_+(A;A)$}\index{$\C^\bullet_-(A;A)$}
\begin{proof}
Straightforward.
\end{proof}

\noindent
These sub-complexes are called 
\emph{the even}\index{H{\nm ochschild}!complex!even part of} and 
\emph{the odd}\index{H{\nm ochschild}!complex!odd part of}
part of $\C^\bullet(A;A)$ respectively.
Sometimes we refer to elements of $\C^\bullet_+(A;A)$ as of 
\emph{even parity}\index{parity of a cochain} and to to elements of 
$\C^\bullet_-(A;A)$ as of \emph{odd parity}.

The second operator on $\C^\bullet(A;A)$ 
\emph{the alternation}\index{alternation} $\Alt$ is defined as
\begin{equation}\label{eq:DefOfAlt}
\Alt\xi=\frac{1}{p!}
\sum_{\si\in\mathfrak{S}_p}(-1)^\si\xi\comp\cP_{\si}
\end{equation}\index{$\Alt$}
where $\xi\in\C^p(A;A)$, $\mathfrak{S}_p$ is the permutation group of
order $p$, $(-1)^\si$ is equal to $+1$ if the permutation $\si$ is 
even and $-1$ if $\si$ is odd,
$\cP$ is the representation of $\mathfrak{S}_p$ in the vector space 
$A^{\ot p}$ defined by 
$\cP_\si(a_1\ot\dots\ot a_p)=
a_{\si^{-1}(1)}\ot\dots\ot a_{\si^{-1}(p)}$.

An element $w\in\C^1(A;A)$ is called 
\emph{derivation}\index{derivation} if $w(ab)=aw(b)+bw(a)$ for all 
$a,b\in A$ ({\nm Leibniz} rule).
Any cocycle in $\C^1(A;A)$ is a derivation of $A$ and vice versa.
\begin{defn}\label{def:Polyderivation}
An element $w\in\C^p(A;A)$ is called 
\emph{$p$-derivation}\index{$p$-derivation} if it obeys the
{\nm Leibniz} rule\index{L{\nm eibniz} rule} for each variable, i.e. 
$$
w(a_1,\dots,a^{'}_ia^{''}_i,\dots,a_p)=
a^{'}_iw(a_1,\dots,a^{''}_i,\dots,a_p)+
a^{''}_iw(a_1,\dots,a^{'}_i,\dots,a_p)
$$ 
for all $i=1,\dots,p$.
\end{defn}
\begin{prop}
Any $p$-derivation is a Hochschild cocycle.
\end{prop}
\begin{proof}
Straightforward computation.
\end{proof}
\begin{prop}\label{th:AltOfHochCoboundIsZero}
Let $A$ be a commutative algebra, $\xi\in\C^p(A;A)$.
Then $\Alt\left(\di\xi\right)=0$.
\end{prop}
\begin{proof}
Using definition (\ref{eq:DefHochschCobound}) and the additivity of
operator $\Alt$, one has:
\begin{eqnarray*}
\Alt\left(\di\xi\right)&=&
\Alt\Big(\mu\comp(\id\ot\xi)+
(-1)^{n+1}\mu\comp(\xi\ot\id)\Big)+ \\
&+&\sum_{k=1}^n(-1)^k\Alt
\Big(\xi\comp(\id^{\ot k-1}\ot\mu\ot\id^{\ot p-k})\Big).
\end{eqnarray*}
The symmetry of $\mu$ implies:
$$
\xi\comp(\id^{\ot k-1}\ot\mu\ot\id^{\ot p-k})=
\xi\comp(\id^{\ot k-1}\ot\mu\ot\id^{\ot p-k})\comp
\cP_{(k,k+1)}
$$
where $(k,k+1)$ is the permutation switching the $k$-th and 
$(k+1)$-th elements and leaving all the others unchanged.
Then 
\begin{eqnarray*}
\Alt\Big(\xi\comp(\id^{\ot k-1}\ot\mu\ot\id^{\ot p-k})\Big)&=&
\Alt\Big(\xi\comp(\id^{\ot k-1}\ot\mu\ot\id^{\ot p-k})\comp
\cP_{(k,k+1)}\Big)=\\
&=&-\Alt\Big(\xi\comp(\id^{\ot k-1}\ot\mu\ot\id^{\ot p-k})\Big)
\end{eqnarray*}
so $\Alt\Big(\xi\comp(\id^{\ot k-1}\ot\mu\ot\id^{\ot p-k})\Big)=0$
for every $k=1,\dots,n$.

It is left to consider the terms
$\Alt\Big(\mu\comp(\id\ot\xi)+(-1)^{p+1}\mu\comp(\xi\ot\id)\Big)$.
Using the symmetry of $\mu$, one has
\begin{eqnarray*}
\mu\comp(\xi\ot\id)(a_1,\dots,a_{p+1})&=&
\mu\comp(\id\ot\xi)(a_{p+1},a_1,\dots,a_{p}) \\
&=&\mu\comp(\id\ot\xi)\comp\cP_\pi(a_1,\dots,a_{p+1})
\end{eqnarray*}
where $\pi=(1\ 2\dots p+1)$, the cyclic permutation.
Since $(-1)^\pi=(-1)^p$, one has
$\Alt\left(\mu\comp(\id\ot\xi)+
(-1)^{p+1}\mu\comp(\xi\ot\id)\right)=0$.
\end{proof}
\begin{prop}\label{th:SkewHochCoboundAndLeibnizRule}
Any skew symmetric Hochschild cocycle is a polyderivation.
\end{prop}
\begin{proof}
Let $\xi\in\C^p(A;A)$ be skew symmetric and suppose $\di\xi=0$.
Since $\xi$ is skew symmetric, it suffices to prove that
$\xi(a_1a_2,a_3,\ldots,a_{p+1})=$\linebreak
$=a_1\xi(a_2,a_3,\ldots,a_{p+1})+a_2\xi(a_1,a_3,\ldots,a_{p+1})$ for
any $a_1,\dots,a_{p+1}\in A$. 
Direct computations show that
$$
\begin{array}{l}
\di\xi(a_1,a_2,a_3,\ldots,a_{p+1})+
\di\xi(a_2,a_1,a_3,\ldots,a_{p+1})=\\
\ \\
=a_1\xi(a_2,a_3,\ldots,a_{p+1})+a_2\xi(a_1,a_3,\ldots,a_{p+1})-
2\xi(a_1a_2,a_3,\ldots,a_{p+1})+\\
\ \\
+\xi(a_1,a_2a_3,\ldots,a_{p+1})+\xi(a_2,a_1a_3,\ldots,a_{p+1})
\end{array}
$$
and
$$
\begin{array}{l}
\Big(\di\xi(a_1,a_2,a_{p+1},a_3,\ldots,a_p)+
\di\xi(a_2,a_1,a_{p+1},a_3,\ldots,a_p)\Big)+\\
\ \\
+\Big(\di\xi(a_{p+1},a_1,a_2,a_3,\ldots,a_p)+
\di\xi(a_{p+1},a_2,a_1,a_3,\ldots,a_p)\Big)=\\
\ \\
=(-1)^p\Big(
a_1\xi(a_2,a_3,\ldots,a_{p+1})+a_2\xi(a_1,a_3,\ldots,a_{p+1})-\\
\ \\
-\xi(a_1,a_2a_3,\ldots,a_{p+1})-
\xi(a_2,a_1a_3,\ldots,a_{p+1})\Big).
\end{array}
$$
Therefore
$$
\begin{array}{l}
0=\di\xi(a_1,a_2,a_3,\ldots,a_{p+1})+
\di\xi(a_2,a_1,a_3,\ldots,a_{p+1})+\\
\ \\
+(-1)^p\Big(\di\xi(a_1,a_2,a_{p+1},a_3,\ldots,a_p)+
\di\xi(a_2,a_1,a_{p+1},a_3,\ldots,a_p)+\\
\ \\
+\di\xi(a_{p+1},a_1,a_2,a_3,\ldots,a_p)+
\di\xi(a_{p+1},a_2,a_1,a_3,\ldots,a_p)\Big)=\\
\ \\
=2\Big(
a_1\xi(a_2,a_3,\ldots,a_{p+1})+a_2\xi(a_1,a_3,\ldots,a_{p+1})-
\xi(a_1a_2,a_3,\ldots,a_{p+1})\Big).
\end{array}
$$
\end{proof}

Later (see Section~\ref{subsect:ExistenceOfFHbPhiHb}) we shall also 
use the {\nm Hochschild} complex for a non-commutative algebra with 
coefficients in two-sided modules.

\subsection{The local Hochschild complex}

Recall that 
\emph{the ring of differential operators on $A=\C^\infty(M)$}
\index{ring of differential operators} is the 
associative algebra generated over $A$ by all derivations $A\to A$
and all multiplication operators $L_a:A\to A,\ b\mapsto ab$.
A map $A^{\ot p}\to A$ is called 
\emph{a $p$-differential operator on $A$}\index{$p$-differential
operators} if it is a differential operator $A\to A$ with respect to 
every its variable.
Denote by $\wtl{\C}^p(A;A)$ the $\kk$ space of $p$-differential
operators on $A$.\index{$\wtl{\C}^p(A;A)$}
\index{$\wtl{\C}^\bullet(A;A)$}
The direct sum of these spaces forms the sub-complex 
$(\wtl\C^\bullet(A;A),\di\,)$ of the complex $(\C^\bullet(A;A),\di\,)$,
it is called 
\emph{the local Hochschild complex of $A$}
\index{H{\nm ochschild}!complex!local}.
We need the local {\nm Hochschild} complex because we deal actually 
with the whole sheaf of function algebras on $M$ rather than with
the algebra $A$ of global sections.
Clearly, all propositions of Section~\ref{subsect:HochCompl} are
valid for the local {\nm Hoch\-schild} complex.
In particular, this complex can be decomposed into the sum of even
and odd parts,
$\wtl\C^\bullet(A;A)=\wtl\C_+^\bullet(A;A)\dsum\wtl\C_-^\bullet(A;A)$
where 
$\wtl\C_\pm^\bullet(A;A)=\C_\pm^\bullet(A;A)\cap\wtl\C^\bullet(A;A)$.

Denote by $\La_p(M)$ the space of all $p$-vector fields on $M$, i.e.
all $p$-derivations $A^{\ot p}\to A$.\index{$\La_p(M)$}
\begin{thm}\label{th:KostHochRosenbKonc}
Let $A=\C^\infty(M)$, $\tht\in\Z^p(\wtl{\C}^\bullet(A;A),\di\,)$.
Then\nl
\mbox{\hspace{.15in}}\emph{(I)}
$\Alt\tht$ is a skew symmetric $p$-derivation on $A$;\nl
\mbox{\hspace{.13in}}\emph{(II)}
The difference $\tht-\Alt\tht$ is a Hochschild coboundary,
i.e. there exists $\xi\in\wtl\C^{p-1}(A;A)$ such that 
$\tht-\Alt\tht=\di\xi$;\nl
\mbox{\hspace{.11in}}\emph{(III)}
The map 
$\Alt:\Z^p(\wtl\C^\bullet(A;A),\di\,)\to\La_p(M)$ induces a
vector space isomorphism 
$\Ha^p(\wtl\C^\bullet(A;A),\di\,)\to\La_p(M)$.
\end{thm}
\begin{proof}
The algebraic version of  (I) and (II) has been proven by 
{\rf  G.~Hochschild}, {\rf  B.~Kostant} and {\rf A.~Rosenberg}
\cite{HochschildKostantRosenberg:DiffFormsOnAffAlgs}.
A proof for $A=\C^\infty(M)$ see in 
\cite{KontsevichM:DefQztnPoissMan1}, 4.6.1.1.
Claim (III) has been proven by {\rf J.~Vey}
\cite{VeyJ:DeformCrochPoiss}, see also 
\cite{LichnerowiczA:AstProduits}.
\end{proof}

\begin{cor}
Let $\tht\in\Z^p(\wtl{\C}^\bullet(A;A),\di\,)$, then $\Alt\tht=0$ implies
that $\tht$ is a Hochschild coboundary.
\end{cor}

\subsection{Infinitesimal of a deformation}
\label{subsect:DefmInfinit}

It was shown above that the associativity of a deformation 
$\mh=\mu+\hb\mu_1+\cdots$ forces the infinitesimal $\mu_1$ to be 
a {\nm Hochschild} cocycle.
It can be decomposed as 
$\mu_1=\mu_1^{+}+\mu_1^{-}\in\C^2_+(A;A)\dsum\C^2_-(A;A)$
(see Proposition~\ref{th:TauIsComplexMorphism}) with
$\mu_1^{-}=\Alt\mu_1$ and $\mu_1^{+}=\mu_1-\mu_1^{-}$.
\begin{prop}
If $\di\mu_1=0$ then both $\di\mu_1^{+}=0$ and $\di\mu_1^{-}=0$.
\end{prop}
\begin{proof}
Clearly, $\di\mu_1=0$ implies $\di\ (\mu_1\comp\cP_{(12)})=0$.
On the other hand,
$\mu_1^{+}=\frac{1}{2}(\mu_1+\mu_1\comp\cP_{(12)})$
and $\mu_1^{-}=\frac{1}{2}(\mu_1-\mu_1\comp\cP_{(12)})$.
\end{proof}

The symmetric part, $\mu_1^{+}$, is responsible for commutative 
deformations of $(A,\mu)$.
For $A=\C^\infty(M)$ any deformation with $\mu_1^{+}+\mu_1^{-}$ as 
the infinitesimal is equivalent to a deformation with infinitesimal 
$\mu_1^{-}$ as the infinitesimal.
Therefore one can assume from the beginning that $\mu_1^{+}=0$
i.e. that the infinitesimal $\mu_1$ is skew symmetric.
By Proposition~\ref{th:SkewHochCoboundAndLeibnizRule}, 
$\di\mu_1=0$ implies that $\mu_1$ is a bivector field.

To simplify our considerations, throughout the text we put the 
following restriction on any deformation $\mh$  
\cite{LichnerowiczA:NATO}.
\begin{defn}\label{def:ParityConvention}
We say that a deformation $\mser$  obeys 
\emph{the Parity Convention}\index{parity convention}
if $\mu_{2k}\in\wtl\C_-^2(A;A)$ and $\mu_{2k+1}\in\wtl\C_+^2(A;A)$.
\end{defn}
\begin{lm}
Let $(\Ah,\mh)$ be a deformation of commutative algebra $A$ obeying
the Parity Convention.
Then its infinitesimal $\mu_1$ satisfies the Jacobi identity.
\end{lm}
\begin{proof}
For $n=2$ equation (\ref{eq:CartMaur}) takes the form
\begin{equation}\label{eq:CartMaurNTwo}
\mu_1\comp(\mu_1\ot\id)-\mu_1\comp(\id\ot\mu_1)=\di\mu_2.
\end{equation}
The right hand side of (\ref{eq:CartMaur}) is a $3$-coboundary.
Straightforward computation with the use of commutativity of $A$ and
symmetry of $\mu_2$ proves that $\di\mu_2(a,b,c)+\mathrm{Cycl}=0$
for any $a,b,c\in A$, where Cycl denotes the terms obtained by 
taking of all the cyclic permutations of $a,b,c$.
For the left hand side of (\ref{eq:CartMaur}), direct computations
with the use of the skew symmetry of $\mu_1$ show that the expression
$\mu_1(\mu_1(a,b),c)-\mu_1(a,\mu_1(b,c))+\mathrm{Cycl}$ is equal
to the left hand side of the {\nm Jacobi} identity.
\end{proof}
\begin{defn}\label{def:QzationOfPoissManif}
\emph{A quantization}\index{quantization!of a {\nm Poisson} manifold}
of a commutative algebra $(A,\mu)$ with a given  element 
$b\in\C_-^2(A;A)$ is a non-commutative deformation of $(A,\mu)$ with 
$b$ as the infinitesimal.
\end{defn}
The  arguments of this section show us that the infinitesimal of any 
quantization is of necessity a {\nm Poisson} bracket.

\newpage
\ \ 
\newpage
\ \ 
\section{Equivalence of monoidal categories\\ related to Drinfeld 
algebras}
\label{sect:MonCategDrAlgs}

In this section we reduce the problem of associative quantization 
invariant under the {\nm Etingof---Kazhdan} quantum group $\Ughr$ to 
the problem of non-associative $G$ invariant quantization.
To do this, we consider two monoidal categories associated to 
the problem of invariant quantization of function algebras, and,
for any {\nm Belavin---Drinfeld} r-matrix, we construct the 
equivalence for these two categories.

Throughout this section, $R$ denotes a commutative associative
ring with unity of characteristic $0$ (we shall be interested in 
cases $R=\kk$ and $R=\kk\dhb$).
All algebras and modules are defined over $R$.

\subsection{Monoidal categories and Drinfeld algebras}

The concept of monoidal category was introduced by {\rf J.~B\'enabou}
\cite{BenabouJ:MonoidCateg} and {\rf S.~Mac~Lane} 
\cite{MacLaneS:NatAssocCommut}, see \cite{MacLaneS:Categories},
Chapter~VII.
\begin{defn}
\emph{A monoidal category}\index{monoidal!category} is a $6$-tuple 
$(\mathfrak{C},\ot,\un,\ba,\bl,\br)$ where $\mathfrak{C}$ is a 
category, $\ot$ is a functor 
$\mathfrak{C}\times\mathfrak{C}\to\mathfrak{C}$,
$\ba$ is a natural isomorphism of functors 
$(U,V,W)\mapsto(U\ot V)\ot W$ and $(U,V,W)\mapsto U\ot(V\ot W)$,
$\bl$ is a natural isomorphism of the functor $V\mapsto \un\ot V$ and 
the identity functor $\Id$,
$\br$ is a natural isomorphism of functors $V\mapsto V\ot\un$ and 
$\Id$,
and the following two diagrams are commutative for any 
$U,V,W,S\in\mathrm{Obj}(\mathfrak{C})$:

\begin{picture}(300,120)(-25,-20)
\put(15,70){\makebox(0,0){$((U\ot V)\ot W)\ot S$}}
\put(68,70){\vector(1,0){30}}
\put(83,80){\makebox(0,0){$\ba$}}
\put(155,70){\makebox(0,0){$(U\ot V)\ot (W\ot S)$}}
\put(207,70){\vector(1,0){30}}
\put(220,80){\makebox(0,0){$\ba$}}
\put(292,70){\makebox(0,0){$U\ot (V\ot (W\ot S))$}}

\put(17,55){\vector(0,-1){40}}
\put(-5,35){\makebox(0,0){$\ba\ot\id$}}
\put(288,15){\vector(0,1){40}}
\put(309,35){\makebox(0,0){$\id\ot\ba$}}

\put(292,1){\makebox(0,0){$U\ot((V\ot W)\ot S)$}}
\put(83,1){\vector(1,0){140}}
\put(145,10){\makebox(0,0){$\ba$}}
\put(17,1){\makebox(0,0){$(U\ot(V\ot W))\ot S$}}
\end{picture}

\begin{picture}(200,120)(10,0)
\put(100,70){\makebox{$U\ot(\un\ot V)$}}
\put(170,70){\vector(1,0){40}}
\put(215,70){\makebox{$(U\ot\un)\ot V$}}
\put(138,60){\vector(1,-1){30}}
\put(243,60){\vector(-1,-1){30}}
\put(173,10){\makebox{$U\ot V$}}
\end{picture}

\noindent
where $\ba$ over the left upper arrow means $\ba_{U\ot V,W,S}$ and so
on.
The isomorphism $\ba$ is called 
\emph{the associativity constraint}\index{associativity constraint}
of $\mathfrak{C}$.
\end{defn}
\begin{defn}
\emph{A monoidal functor}\index{monoidal!functor} from a monoidal 
category $(\mathfrak{C},\ot,\ba)$ to a monoidal category 
$(\mathfrak{D},\boxtimes,\bb)$ is a pair $(L,\bu)$ 
where $L$ is a functor $\mathfrak{C}\to\mathfrak{D}$, 
$\bu$ is a natural transformation from the functor 
$(U,V)\mapsto L(U\ot V)$ to $(U,V)\mapsto L(U)\boxtimes L(V)$ such 
that $\bu\comp L(\bl_\mathfrak{C})=\bl_\mathfrak{D}$,
$\bu\comp L(\br_\mathfrak{C})=\br_\mathfrak{D}$
and the following diagram is commutative:

\begin{picture}(300,120)(0,-20)
\put(15,70){\makebox(0,0){$L((U\ot V)\ot W)$}}
\put(68,70){\vector(1,0){30}}
\put(80,80){\makebox(0,0){$\bu$}}
\put(155,70){\makebox(0,0){$L(U\ot V)\boxtimes L(W)$}}
\put(210,70){\vector(1,0){50}}
\put(235,80){\makebox(0,0){$\bu\boxtimes\id$}}
\put(325,70){\makebox(0,0)
    {$\left(L(U)\boxtimes L(V)\right)\boxtimes L(W)$}}

\put(17,55){\vector(0,-1){40}}
\put(-3,35){\makebox(0,0){$L(\ba)$}}
\put(325,55){\vector(0,-1){40}}
\put(340,35){\makebox(0,0){$\bb$}}

\put(15,1){\makebox(0,0){$L(U\ot(V\ot W))$}}
\put(68,1){\vector(1,0){30}}
\put(80,10){\makebox(0,0){$\bu$}}
\put(155,1){\makebox(0,0){$L(U)\boxtimes L(V\ot W)$}}
\put(210,1){\vector(1,0){50}}
\put(235,10){\makebox(0,0){$\id\boxtimes\bu$}}
\put(325,1){\makebox(0,0)
     {$L(U)\boxtimes\left(L(V)\boxtimes L(W)\right)$}}
\end{picture}

\end{defn}
\begin{eg}
Let $(B,\De,\ep)$ be a bialgebra over\  $\kk$.
Consider the category $\mathfrak{C}_B$ of all modules over $B$
which are of finite dimension as $\kk$ vector spaces.
The category $\mathfrak{C}_B$ possesses a monoidal structure:
for given $B$ modules $M$ and $N$ one can introduce a $B$
module structure on their tensor product $M\ot N$ over $\kk$ by 
putting $b(m\ot n)=\sum b'm\ot b''n$  where $b\in B$,
$m\ot n\in M\ot N$ and $\De(b)=\sum b'\ot b''$.
\end{eg}
{\rf V.~G.~Drinfeld} has generalized the above example in the 
following way.
\begin{defn}\label{def:DrinfeldAlg}
A \emph{Drinfeld algebra}\index{D{\nm rinfeld} algebra} is a 
$6$-tuple $(B,m,\iota,\De,\ep,\Ph)$ where $(B,m,\iota)$ is an 
associative $R$ algebra with a unit map $\iota:R\to B$,
$(B,\De,\ep)$ is a coalgebra with comultiplication 
$\De:B\to B\ot_R B$ and counit $\ep:B\to R$ are algebra morphisms,
$\Ph\in B\ot_RB\ot_RB$ is an invertible element,
and the following conditions are obeyed,
where we use the notation $a\cdot b=m(a,b)$:\nl
\emph{(i)} $(\id\ot\De)\comp\De(b)\cdot\Ph=\Ph\cdot
(\De\ot\id)\comp\De(b)$ for any $b\in B$;\nl
{\em(ii)} $(\id\ot\id\ot\De)(\Ph)\cdot
(\De\ot\id\ot\id)(\Ph)=
(1\ot\Ph)\cdot(\id\ot\De\ot\id)(\Ph)\cdot
(\Ph\ot 1)$;\nl
\emph{(iii)}
$(\ep\ot\id)\comp\De=\id=
(\id\ot\ep)\comp\De$;\nl
\emph{(iv)}
$(\id\ot\ep\ot\id)(\Ph)=1$.\nl
The element $\Ph$ is called 
\emph{the Drinfeld associator of $B$}\index{D{\nm rinfeld}
associator}.
The condition \emph{(ii)} is called the 
\emph{pentagon}\index{pentagon identity} or
\emph{Mac Lane identity}\index{M{\nm ac Lane} identity} for $\Ph$.
\end{defn} 
\begin{eg}
The trivial associator, $\Ph=1\ot 1\ot 1$, obeys all the 
conditions of Definition~\ref{def:DrinfeldAlg}.
In this case $B$ is a usual coassociative bialgebra.
\end{eg}
\begin{rem}
We did not put the antipode (\cite{AbeE:HopfAlgs},Chapter~2,
Subsection~1.2) in the above definition,
because we will not referring to it in any of our results.
Yet in the main example, the QUE \emph{Drinfeld} algebra 
(Definition~\ref{def:QUEDrinfAlgs}) the antipode exists as a
deformation of the standard antipode on $\Ug$.
\end{rem}
{\nm Drinfeld} algebras were introduced in 
\cite{DrinfeldW:qHopfAlgs} under the name 
\emph{quasi Hopf algebras}\index{quasi {\nm Hopf} algebra}.
For a given {\nm Drinfeld} algebra $(B,\Ph)$ over $R$, denote by 
$\mathfrak{D}_B$ the category of $B$ modules which are free $R$
modules of finite rank.
Note that $\mathfrak{D}_B$ is a sub-category of the category
$\mathfrak{D}_R$.
In particular, one can consider in $\mathfrak{D}_B$ tensor 
products over $R$.
\begin{prop}
The category $\mathfrak{D}_B$ possesses a monoidal structure 
with associativity constraint defined by $\Ph$.
\end{prop}
\begin{proof}
Consider the usual tensor product $(U,V)\mapsto U\ot_RV$ over $R$
in $\mathfrak{D}_B$.
One can equip $U\ot_RV$ with a $B$-action by setting 
$b.(u\ot v)=\De(b)\cdot(u\ot v)$ for any
$b\in B,\ \ u\in U,\ \ v\in V$.
Denote this $B$ module by $U\ot_B V$.

Write $\Ph$ as $\sum\Ph^{'}\ot\Ph^{''}\ot\Ph^{'''}$
and for any $u\in U,\ v\in V,\ w\in W$ define
$$
\mathbf{a}_{UVW}\Big((u\ot v)\ot w\Big)=
\sum(\Ph^{'}\cdot u)\ot\Big((\Ph^{''}\cdot v)
\ot(\Ph^{'''}\cdot w)\Big).
$$
One can check that $(\mathfrak{D}_B,\ot_B,\mathbf{a})$ is a monoidal 
category.
\end{proof}

\noindent
The difference between {\nm Drinfeld} algebras and the usual 
bialgebras is that the associativity constraint in the category
of modules over a {\nm Drinfeld} algebra is in general non-trivial.
We shall quantize in a category over a {\nm Drinfeld} algebra,
i.e. we shall look for a $\Ph$-associative deformed multiplication
rather than for a multiplication with the ordinary associativity.
\begin{defn}
Monoidal categories $\mathfrak{C}$ and $\mathfrak{D}$ are called 
\emph{equivalent}\index{monoidal category!equivalence} if there 
exist two monoidal functors, $F:\mathfrak{C}\to\mathfrak{D}$ and 
$G:\mathfrak{D}\to\mathfrak{C}$ such that the both compositions of
$F$ and $G$ are naturally equivalent to the corresponding 
identity functors.
\end{defn}
\begin{thm}\label{th:GaugeTransform}
Let $(B,\De,\Ph)$ be a Drinfeld algebra, $F\in B\ot_R B$ an 
invertible element.
Define $\wtl\De$ and $\wtl\Ph$ by 
$\wtl\De(b)=F\cdot\De(b)\cdot F^{-1}$ for all $b\in D$ and 
$\wtl\Ph=
(1\ot F)\cdot(\id\ot\De)(F)\cdot\Ph\cdot
(\De\ot\id)(F^{-1})\cdot(F\ot 1)^{-1}$.
Then $(B,m,\iota,\wtl\De,\ep,\wtl\Ph)$ is also a {\nm Drinfeld} 
algebra and the categories $\mathfrak{D}_B$ and 
$\mathfrak{D}_{\wtl{B}}$ are monoidally equivalent.
\end{thm}
\begin{proof}
See \cite{DrinfeldW:qHopfAlgs}.
\end{proof}

\noindent
Note that gauge transform is called in \cite{DrinfeldW:qHopfAlgs} 
\emph{the twisting by $F$}\index{twisting}.

\begin{defn}\label{def:GaugeTransform}
The Drinfeld algebra $(B,m,\iota,\wtl\De,\ep,\wtl\Ph)$
of Theorem~\ref{th:GaugeTransform} is called 
\emph{a gauge transformation}\index{gauge transformation} of 
$(B,m,\iota,\De,\ep,\Ph)$  \emph{determined by $F$}.
Two Drinfeld algebras are called 
\emph{gauge equivalent}\index{gauge equivalence} if there exists a 
gauge transformation from one to the other.
\end{defn}
Theorem~\ref{th:GaugeTransform} states that if two {\nm Drinfeld}
algebras are gauge equivalent then the corresponding monoidal
categories are equivalent.
The equivalence is given by the pair $(\Id,F)$.

\subsection{Monoids in a category and invariant multiplications}

The monoidal structure in a category makes it possible to construct 
algebra-like objects in it.
In the category $\mathfrak{D}_B$, the corresponding multiplication 
is called $B$ invariant.
\begin{defn}\label{def:Monoid}
\emph{A monoid}\index{monoid} in a monoidal category
$(\mathfrak{C},\ot,\un,\ba)$ is a triple $(A,\mu,\eta)$ 
where $A\in \mathrm{Obj}({\mathfrak C})$, 
$\mu:A\ot A\to A$, and $\eta:\un\to A$ are arrows of $\mathfrak{C}$,
and the following three diagrams are commutative:

\begin{picture}(300,120)(-100,-17)
\put(15,70){\makebox(0,0){$(A\ot A)\ot A$}}
\put(115,70){\makebox(0,0){$A\ot A$}}
\put(15,10){\makebox(0,0){$A\ot(A\ot A)$}}
\put(115,10){\makebox(0,0){$A\ot A$}}
\put(53,70){\vector(1,0){40}}
\put(51,10){\vector(1,0){45}}
\put(70,80){\makebox(0,0){$\mu\ot\id_A$}}
\put(70,0){\makebox(0,0){$\id_A\ot\mu$}}
\put(167,40){\makebox(0,0){$A$}}
\put(15,55){\vector(0,-1){30}}
\put(135,65){\vector(1,-1){22}}
\put(135,11){\vector(1,1){22}}
\put(5,40){\makebox(0,0){$\ba$}}
\put(150,65){\makebox(0,0){$\mu$}}
\put(150,15){\makebox(0,0){$\mu$}}
\end{picture}

\begin{picture}(300,120)(140,-10)
\put(200,70){\makebox(0,0){$\un\ot A$}}
\put(265,70){\makebox(0,0){$A\ot A$}}
\put(217,70){\vector(1,0){30}}
\put(222,78){\makebox{$\eta\ot\id$}}
\put(200,10){\makebox(0,0){$A$}}
\put(198,58){\vector(0,-1){35}}
\put(263,30){\makebox(0,0){$\mu$}}
\put(265,58){\vector(-1,-1){50}}

\put(370,70){\makebox(0,0){$A\ot\un$}}
\put(435,70){\makebox(0,0){$A\ot A$}}
\put(387,70){\vector(1,0){30}}
\put(393,78){\makebox{$\id\ot\eta$}}
\put(370,10){\makebox(0,0){$A$}}
\put(368,58){\vector(0,-1){35}}
\put(433,30){\makebox(0,0){$\mu$}}
\put(435,58){\vector(-1,-1){50}}
\end{picture}
\end{defn}
\begin{eg}
Let $B$ be a bialgebra.
Then a monoid $(A,\mu)$ in the category ${\mathfrak{D}}_B$ is a usual 
associative $B$ algebra with invariance property 
$x.\mu(a,b)=\mu(x^{'}a,x^{''}b)$, where $\De(x)=\sum x^{'}\ot x^{''}$.
\end{eg}
\begin{eg}\label{eg:PhAssoc}
Let $(B,\De,\Ph)$ be a {\nm Drinfeld} algebra where 
$\Ph=\sum_\Ph\Ph'\ot\Ph''\ot\Ph'''$.
Then a monoid $(A,\mu)$ in ${\mathfrak{D}}_B$ is a $B$ algebra 
with a multiplication $a\cdot b=\mu(a\ot b)$ obeying
$$
(a\cdot b)\cdot c=\sum_\Ph\Ph'a\cdot(\Ph''b\cdot\Ph'''c)
$$
rather than the usual associativity.
Also this can be written as 
$\mu\comp(\mu\ot\id)=\mu\comp(\id\ot\mu)\comp\Ph$
where $\Ph$ is considered as a left multiplication operator.
We call the multiplication law of such $A$ \emph{a $B$ invariant
multiplication}\index{multiplication!invariant} or 
\emph{$\Ph$-associative multiplication}.\index{$\Ph$-associativity}
\end{eg}
\begin{defn}\label{def:GInvarMult}
Let $\g$ be a Lie algebra,
$\Ug$ the universal enveloping algebra of $\g$,
$\Ph\in\Ug^{\ot 3}$ an element satisfying
the conditions of Definition~\ref{def:DrinfeldAlg}.
Any multiplication in monoid in the category 
${\mathfrak{D}}_{(\Ug,\Ph)}$ we call 
\emph{ $\g$ invariant}\index{multiplication!$\g$ invariant}.
\end{defn}
\begin{eg}
Consider $\Ug$ as a {\nm Drinfeld} algebra with the trivial
associator.
The algebra $(A,\mu)$ of $\C^\infty$ functions on $M$ is $\Ug$ 
invariant if and only if $x\comp\mu=\mu\comp\De(x)$ for any 
$x\in\Ug$ considering as a differential operator on $A$.
Since the comultiplication $\De$ is an algebra $(\Ug,m)$
morphism,
it suffices to check the above identity on elements of any set of 
generators of $\Ug$.
For instance, one can take elements of the Lie algebra $\g$ itself:
$x(ab)=\mu\left(\De(x).(a\ot b)\right)=
\mu\left((x\ot 1+1\ot x).(a\ot b)\right)=
x(a)b+ax(b)$.
Thus the multiplication $\mu$ is $\g$ invariant if and only if
elements of $\g$ act on $A$ by derivations.
\end{eg}

\subsection{Belavin---Drinfeld r-matrices and QUE Drinfeld
algebras}

Let $\g$ be a complex simple {\nm Lie} algebra.
It is well known that the $\g$ invariant elements of $\bw^3\g$ form 
a one dimensional subspace.
A base is given by the element $\ph$\index{$\ph$} representing 
$x\ot y\ot z\mapsto ([x,y],z)$ where $(\ \cdot\ ,\ \cdot\ )$ is the 
{\nm Killing} form in $\g$.
\begin{defn}\label{def:mCYBE}
The equation $\dl r,r\dr=\ph$ is called
\emph{the modified classical Yang---Baxter equation}.
\index{Y{\nm ang--Baxter} classical equation!modified}
\index{modified CYBE}
Any solution of it is called \emph{a Belavin---Drinfeld r-matrix on 
$\g$}\index{B{\nm elavin--Drinfeld}!r-matrix}
\index{r-matrix!B{\nm elavin--Drinfeld}}.
\end{defn}
\begin{eg}
Let $\Wurz^+$ be a system of positive roots for $\g$,
$X_\al$ corresponding root vectors satisfying $(X_\al,X_{-\al})=1$.
Then the element of $\bw^2\g$ defined by
$$
r=\sum_{\al\in\Wurz^+}X_\al\w X_{-\al}
$$
is a Belavin---Drinfeld r-matrix.
This element is called \emph{the Drinfeld---Jimbo r-matrix}
\index{D{\nm rinfeld--Jimbo} r-matrix}
\index{r-matrix!{\nm Drinfeld--Jimbo}} and is
the most important example of a {\nm Belavin--Drinfeld} r-matrix.
\end{eg}
All solutions for the modified classical {\nm Yang---Baxter} equation 
were found by {\rf A. A. Belavin} and {\rf V. G. Drinfeld}
\cite{BelawinDrinfeld:CYBEpreprint}, \cite{BelawinDrinfeld:CYBE}.
The significance of {\nm Belavin---Drinfeld} r-matrices in the
fact \cite{DrinfeldW:CYBEgeommeaning} that any such element 
$r\in\bw^2\g$ determines a {\nm Lie} algebra structure on the
dual vector space $\gd$.
This structure is consistent with the {\nm Lie} bracket on $\g$, i.e.
$r$ is a {\nm Chevalley--Eilenberg} 2-coboundary.
The pair $(\g,r)$ is called \emph{a Belavin---Drinfeld Lie 
bialgebra}\index{B{\nm elavin--Drinfeld}!{\nm Lie} bialgebra}.
A {\nm Belavin---Drinfeld} $r$-matrix is a natural initial term
for a coalgebra deformation,
as a {\nm Poisson} structure is a natural initial term for an
algebra deformation.
\begin{rem}\label{nb:HTwoGvrnsIsomClss}
It follows from the M.~Gerstenhaber's obstruction theory
\cite{GerstenhaberM:DeformsRingsAlgs1} that since 
$\Ha^2(\Ug;\Ug)=0$ for $\g$ semisimple, there are no 
non-trivial deformations of the multiplication in the universal 
enveloping algebra $\Ug$.
\end{rem}
\begin{defn}\label{def:QUEDrinfAlgs}
A deformation $(\Ug\dhb,m,\iota,\Db,\ep,\Pb)$ of the Drinfeld algebra
$(\Ug,m,\iota,\De,\ep,1\ot 1\ot 1)$ is called
\emph{a quantized universal enveloping  (QUE) Drinfeld algebra}.
\index{QUE {\nm Drinfeld} algebra}
\emph{A quantization of a Belavin---Drinfeld Lie bialgebra
$(\g,r)$} is a QUE Drinfeld algebra with the deformed
comultiplication of the form $\Db(x)=\De(x)+\hb[\De(x),r]+\cdots$.
\index{quantization!of {\nm Belavin--Drinfeld Lie} bialgebra}
\end{defn}
Note that we have a normalized form in which the multiplication is
undeformed, which is always possible by 
Remark~\ref{nb:HTwoGvrnsIsomClss}.

\subsection{Uniqueness of Drinfeld quantization for $\g$ simple\\
and existence of $\Fb$ for any Belavin--Drinfeld infinitesimal}
\label{subsect:ExistenceOfFHbPhiHb}

For a given {\nm Lie} algebra $\g$,
denote by $\Ha_{CE}^\bullet\Big(\g;\Ug^{\ot 2}\Big)$ its 
{\nm Chevalley--Eilen\-berg} cohomology module and by 
$\Ha_{Hoch}^\bullet\Big(\Ug;\Ug^{\ot 2}\Big)$ its {\nm Hochschild} 
cohomology module. 
To construct the {\nm Chevalley--Eilenberg} complex,
one uses the action $(X,v)\mapsto [\De(X),v]=\De(X)v-v\De(X)$ of
$\g$ on $\Ug\ot\Ug$,
while in the {\nm Hochschild} complex the right and left
actions determined by the standard comultiplication $\De$ on $\Ug$
are utilized.
\begin{lm}\label{th:CEOneIsEqualHochOne}
Let $\g$ be a finite dimensional Lie algebra.
Then $\Ha^1_{CE}\Big(\g;\Ug^{\ot 2}\Big)=0$ implies 
$\Ha^1_{Hoch}\Big(\Ug;\Ug^{\ot 2}\Big)=0$.
\end{lm}
\begin{proof}
The restriction $\wtl{\xi}$ to $\g$ of a {\nm Hochschild} $1$-cocycle 
$\xi:\Ug\to\Ug\ot\Ug$ is a {\nm Chevalley--Eilenberg} cocycle.
Indeed, $\di_H\xi=0$ implies  
\begin{equation}\label{eq:HochOneCocycle}
\xi(XY)=\De(X)\xi(Y)+\xi(X)\De(Y) \mbox{ for any } X,Y\in\g.
\end{equation}
The relation $XY-YX-[X,Y]=0$ in $\Ug$ 
(where $[X,Y]$ is the {\nm Lie} bracket in $\g$) implies 
$\xi([X,Y])=\xi(XY)-\xi(YX)$.
Substituting  this expression with $(\ref{eq:HochOneCocycle})$,
one has 
$$
\xi([X,Y])=\De(X)\xi(Y)-\xi(Y)\De(X)+\xi(X)\De(Y)-\De(Y)\xi(X)
$$ 
for any $X,Y\in\g$, that is, $\xi$ is a {\nm Chevalley--Eilenberg} 
cocycle.

The above considerations together with the condition
$\Ha^1_{CE}\Big(\g;\Ug\ot\Ug)\Big)=0$ imply that there exists 
$\eta\in\Ug\ot\Ug$ such that $\di_{CE}\eta=\wtl{\xi}$ \  i.e. such that
$\De(X)\eta-\eta\De(X)=\wtl{\xi}(X)$ for any $X\in\g$.
The element $\eta$ can be considered also as a {\nm Hochschild}
$0$-cochain with $\di_H\eta(x)=\De(x)\eta-\eta\De(x)$.
Therefore $\di_H\eta(X)=\di_{CE}\eta(X)$, 
and also $\di_H\eta(X)=\xi(X)$ for any $X\in\g$.

Since both $\xi$ and $\di_H\eta$ are {\nm Hochschild} cocycles,
they obey the condition $(\ref{eq:HochOneCocycle})$.
Combining $(\ref{eq:HochOneCocycle})$ with the fact that $\g$ 
generates the associative algebra $\Ug$,
one concludes that $\wtl\xi$ can be extended uniquely to a
{\nm Hochschild} cocycle,
and that $\di_H\eta(x)=\xi(x)$ for all $x\in\Ug$.
Thus any cocycle in $\C^1_{Hoch}\Big(\Ug;\Ug^{\ot 2}\Big)$ can be 
resolved.
\end{proof}

\begin{lm}\label{th:AbstractFExists}
Let $\de:A\to B$ is a homomorphism of two algebras,
$\de_\hb:A\dhb\to B\dhb$ a deformation of $\de=\de_0$.
Consider the structure of two-sided $A$ module on $B$ 
determined by $\de$ and suppose $\Ha^1_{Hoch}(A;B)=0$.
Then there exists an invertible element $\Fb\in B\dhb$ such
that $\de(x)=\Fb\cdot\de_\hb(x)\cdot\Fb^{-1}$ for any $x\in A\dhb$.
\end{lm}
\begin{proof}
Obviously, the element $1+\hb^kf_k$ is invertible in $B\dhb$ for 
any $f_k\in B$.
Put $\de^{(1)}=\de_\hb$ and prove that for any $k=1,2,\dots$ one can
find  an element $f_k\in B$ such that each algebra homomorphism 
$A\dhb\to B\dhb$ defined by 
$\de^{(k+1)}(a)=(1+\hb f_k)\de^{(k)}(a)(1+\hb f_k)^{-1}$
is of the form $\de+\hb^{k+1}\de_{k+1}^{''}+\cdots$.
Indeed, gathering all the terms with $\hb^k$  the equation
$(1+\hb^k f_k)(\de(a)+\hb^k\de^{'}_k+\hb^{k+1}\de_{k+1}+\cdots)
=(\de(a)+\hb^{k+1}\de_{k+1}^{''}+\cdots)(1+\hb^k f_k)$, one obtains 
the equation $f_k\de(a)-\de(a)f_k=-\de_k^{'}(a)$ for $f_k$.
One checks readily that the left hand side of this equation is
equal to $-\di f_k(a)$ where $\di$ is the {\nm Hochschild}
coboundary operator.
Thus the equation for $f_k$ is of the form $\di f_k=\de_k^{'}$.
Since $\Ha^1_{Hoch}(A;B)=0$, the equation has a solution if 
$\de_k^{'}$ is a {\nm Hochschild} cocycle,
which is the case since $\de^{(k)}=\de+\hb^k\de_k^{'}+\cdots$ is 
multiplicative, i.e. $\de^{(k)}(ab)=\de^{(k)}(a)\de^{(k)}(b)$
for any $a,b\in A$.
For $\de_k^{'}$, this gives the property
$\de_k^{'}(ab)=\de(a)\de_k^{'}(b)+\de_k^{'}(a)\de(b)$ which 
exactly means that $\de_k^{'}$ is a {\nm Hochschild} cocycle.
Finally, denote by $\Fb=\Pi(1\ot 1+\hb^if_i)$, the infinite 
product on the left.
\end{proof}
\begin{prop}\label{th:FExists}
Let $\g$ be a semisimple complex {\nm Lie} algebra,
$\Db:\Ug\dhb\to\Ug^{\ot 2}\dhb$ a deformation of the standard 
comultiplication $\De$ on $\Ug$.
Then there exists an invertible element $\Fb\in\Ug^{\ot 2}\dhb$ such 
that $\Db(x)=\Fb^{-1}\cdot\De(x)\cdot \Fb$ for any $x\in\Ug\dhb$.
\end{prop}\index{$\Fb$}
\begin{proof}
It is well known that for $\g$ semi-simple and for any $\g$ module $V$
one has $\Ha^1_{CE}(\g;V)=0$ (see for instance 
\cite{PostnikovMM:LectGeomV}, Lecture~19, Corollary~1).
Therefore $\Ha_{CE}^1\Big(\g,\Ug^{\ot 2}\Big)=0$,
and by Lemma~\ref{th:CEOneIsEqualHochOne}, 
this implies $\Ha_{Hoch}^1\Big(\Ug,\Ug^{\ot 2}\Big)=0$.
Now one can apply Lemma~\ref{th:AbstractFExists} to the case 
when $B=\Ug\ot\Ug$, $\de=\De$ and $\de_\hb=\Db$.
\end{proof}

{\rf P.~Etingof} and {\rf D.~Kazhdan} have proven 
\cite{EtinhoffKazdan:QuantBialg} that any {\nm Lie} bialgebra can be
quantized.
For a {\nm Belavin--Drinfeld Lie} bialgebra $(\g,r)$ this means that
there exists a bialgebra deformation $\Ughr=(\Ug\dhb,\Db)$ with 
undeformed multiplication (see Remark~\ref{nb:HTwoGvrnsIsomClss}) 
and comultiplication $\Db$ of the form 
$\Db(x)=\De(x)+\hb[\De(x),r]+\cdots$.
We shall call $\Ughr$ \emph{the Etingof--Kazhdan quantum group
corresponding to} $r$.\index{E{\nm tinhof--Kazhdan} quantum group}
\index{quantum group!E{\nm tinhof--Kazhdan}}
In this case, the element $\Fb$ of Proposition~\ref{th:FExists} is of 
the form $\Fb=1\ot 1+\frac{1}{2}\hb r+\cdots$.

\begin{prop}\label{th:BringPhQuadTrmInvForm}
Let $\Pb\in\Ug^{\ot 3}$ be an invertible element $\Pb\equiv 1\mdhs[2]$
satisfying the Pentagon Identity.
Then there exists a gauge transformation making $\Pb$ into
$\wtl\Pb\equiv 1+\hb^2\ph\mdhs[3]$ for some $\ph\in\bw^3\g$.
If $\Pb$ satisfies $\Pb\cdot\De_3(x)\cdot\Pb^{-1}=\De_3(x)$ for the
standard comultiplication $\De$ and for all $x\in\Ug$ 
then the element $\ph$ is $\g$ invariant.
\end{prop}
\begin{proof}
We use the {\nm Cartier} cochain complex\index{C{\nm artier} complex}
$$
\kk\to\Ug\stackrel{\di_C}{\to}\Ug^{\ot 2}\stackrel{\di_C}{\to}
\Ug^{\ot 3}\to\cdots,
$$
$$
\di_C=1\ot\id^{\ot p}+
\sum_{k=1}^p(-1)^k\id^{\ot(k-1)}\ot\De\ot\id^{\ot k}+
(-1)^{k+1}\id^{\ot p}\ot 1.
$$
Let $\Pb=1+\hb^2\psi+\cdots$.
Putting together all the terms of the pentagon identity for 
$(\Db,\Pb)$ which are quadratic in $\hb$, and using 
$\Db(1)=1\ot 1$, one obtains $\di_C\psi=0$.
It is well known (see for example \cite{DrinfeldW:qHopfAlgs},
Proposition~2.2) that the cohomology module for the {\nm Cartier}
complex is equal to $\La(\g)$.
Thus there exist $\ph\in\bw^3\g$ and $f\in\Ug^{\ot 2}$ such that
$\psi=\ph+\di_C f$.
Take $\wtl\Fb=1+\hb^2f$ and 
$\wtl\Pb=(\wtl\Fb\ot 1)\cdot(\Db\ot\id)(\wtl{\Fb})\cdot\Pb
\cdot(\id\ot\Db)(\wtl\Fb^{-1})\cdot(1\ot\wtl\Fb)^{-1}$.
Obviously, $\wtl\Pb\equiv 1+\hb^2\ph\mdhs[3]$.

Suppose now that $\wtl\Pb\cdot\De_3(x)\cdot\wtl\Pb^{-1}=\De_3(x)$ 
for all $x\in\Ug$.
Taking the terms with $\hb^2$ in this equation,
one obtains $\ph\cdot\De_3(x)=\De_3(x)\cdot\ph$ for any $x\in\Ug$
which is an expression for the $\g$ invariance of $\ph$.
\end{proof}

The following theorem follows from results in 
\cite{DrinfeldW:qHopfAlgs}.
\begin{thm}\label{th:QzationBelDrinfMatrix}
Let $r$ be a Belavin--Drinfeld r-matrix,
$\Ughr$ the corresponding Etingof--Kazhdan quantum group,
$m$ and $\De$ the standard multiplication and comultiplication on 
$\Ug$ respectively.
Then there exists an invertible $\g$ invariant element 
$\Pb\in\Ug^{\ot 3}\dhb$ such that\nl
\mbox{\hspace{.2in} \emph{(I)}\hspace{.2in}}
$(\Ug\dhb,m,\De,\Pb)$ is a Drinfeld algebra;\nl
\mbox{\hspace{.17in} \emph{(II)}\hspace{.17in}}
$\Pb\equiv 1+\hb^2\ph\mdhs[3]$ where $0\not=\ph\in(\bw^3\g)^\g$;\nl
\mbox{\hspace{.15in} \emph{(III)}\hspace{.15in}}
There exists an element $\Fb=1+\hb r+\cdots$ which defines a gauge 
equivalence between $\Ughr$ and $(\Ug\dhb,m,\De,\Pb)$.
Therefore the module categories corresponding to $\Ughr$ and 
$(\Ug\dhb,m,\De,\Pb)$ are equivalent as monoidal categories.
\end{thm}\index{$\Ughr$}\index{$\Pb$}\index{$\Fb$}
\begin{proof}
First we prove that there exists an element $\Pb\in\Ug^{\ot 3}\dhb$ 
satisfying $\Pb\equiv 1\mdhs[2]$ and such that $(\Ug\dhb,\De,\Pb)$ is
a {\nm Drinfeld} algebra, gauge equivalent to $\Ughr$.

By Proposition~\ref{th:FExists}, 
there exists $\Fb$ such that $\Fb^{-1}\cdot\De(x)\cdot\Fb=\Db(x)$.
In our case, $\Db$ is coassociative, that is,
\begin{equation}\label{eq:DeHbCoass}
(\Db\ot\id)\comp\Db=(\id\ot\Db)\comp\Db.
\end{equation}
Using multiplicativity of $\De$,
one has for the left hand side of $(\ref{eq:DeHbCoass})$:
$$
\begin{array}{l}
(\Db\ot\id)(\Db(x))
=(\Fb\ot 1)^{-1}\cdot(\De\ot\id)(\Db(x))\cdot(\Fb\ot 1)=\\
\ \\
=(\Fb\ot 1)^{-1}\cdot(\De\ot\id)
(\Fb^{-1}\cdot\De(x)\cdot \Fb)\cdot(\Fb\ot 1)=\\
\ \\
=(\Fb\ot 1)^{-1}\cdot(\De\ot\id)(\Fb^{-1})\cdot\De_3(x)\cdot
(\De\ot\id)(\Fb)\cdot(\Fb\ot 1)
\end{array}
$$
where the notation $\De_3=(\De\ot\id)\comp\De$ is used.
Similarly, the right hand side of $(\ref{eq:DeHbCoass})$ is of the
form $(1\ot \Fb)^{-1}\cdot(\id\ot\De)(\Fb^{-1})\cdot\De_3(x)\cdot
(\id\ot\De)(\Fb)\cdot(1\ot \Fb)$ since 
$(\De\ot\id)\comp\De=(\id\ot\De)\comp\De$.
Put $\Pb=(\id\ot\De)(\Fb)\cdot(1\ot \Fb)
\cdot(\Fb\ot 1)^{-1}\cdot(\De\ot\id)(\Fb^{-1})$.
Clearly, $\Pb$ satisfies $\Pb\cdot\De_3(x)\cdot\Pb^{-1}=\De_3(x)$
for all $x\in\Ug\dhb$.
Considering $(\Ug\dhb,\Db)$ as a {\nm Drinfeld} algebra with
the trivial associator and using Theorem~\ref{th:GaugeTransform},
one concludes that $(\Ug\dhb,\De,\Pb)$ is also a {\nm Drinfeld} 
algebra.

We prove that $\Pb$ has no linear term in $\hb$.
Recall that $\Fb\equiv 1\ot 1+\hb r\mdhs[2]$,
thus $\Fb^{-1}\equiv 1\ot 1-\hb r\mdhs[2]$. 
Using the standard notation ($1\ot r=r^{23}$, $r\ot 1=r^{12}$ and 
$r^{13}=\sum a_i\ot 1\ot b_i$ for $r=\sum a_i\ot b_i$) and the fact 
that $r\in\g\ot\g$, one has $(\id\ot\De)(r)=r^{12}+r^{13}$ and 
$(\De\ot\id)(r)=r^{13}+r^{23}$.
Hence the linear in $\hb$ term of $\Pb$ is equal to
$r^{12}+r^{13}+r^{23}-r^{12}-r^{13}-r^{23}=0$.

Combining this with Proposition~\ref{th:FExists},  
Proposition~\ref{th:BringPhQuadTrmInvForm} 
and Theorem~\ref{th:GaugeTransform}, one completes the proof.
\end{proof}

\begin{rem}\label{nb:MultiplsInEquivCategs}
As we mentioned above,
the gauge transformation $\Fb=1\ot 1+\hb r+\cdots$ from $\Ughr$ to
$(\Ug\dhb,m,\De,\Pb)$ establishes a monoidal equivalence between 
$\mathfrak{C}=\mathfrak{D}_{(\Ug\dhb,m,\De,\Pb)}$ and 
$\mathfrak{C}^{'}=\mathfrak{D}_{(\Ug\dhb,m,\Db,\Id^{\ot 3})}$.
If $(A\dhb,\mh)$ is an algebra in $\mathfrak{C}$ then the
corresponding algebra in $\mathfrak{C}^{'}$ is $(A\dhb,\mh\comp\Fb)$.
In explicit terms, if $\mh=\mu+\hb\mu_1+\cdots$ then 
$\mh\comp\Fb=\mu+\hb(\mu_1+r)+\dots$.
\end{rem}

\newpage
\ \ 
\newpage
\section{Infinitesimals for $\Ughr$ invariant quantizations}
\label{sect:InfinitesForInvarQn}

Let $G$ be a connected simple {\nm Lie} group over $\kk$, $\g$ the
{\nm Lie} algebra of $G$.
In this section we study properties of infinitesimals for $G$
invariant quantization of function algebras on homogeneous manifolds
in the category with {\nm Drinfeld} associator $\Pb$.
We show that, similarly to the case of associative deformation, any 
$\Pb$-associative deformation has as infinitesimal a bracket obeying
the {\nm Leibniz} rule and a weak version of the {\nm Jacobi} 
identity.

Throughout the section, $M$ denotes a homogeneous $G$ manifold, 
thus it is isomorphic to $G/K$ for some closed {\nm Lie} subgroup 
$K$ of $G$.
We fix a \emph{reductive} closed subgroup $K$ and denote by $\k$ its 
{\nm Lie} algebra.
By $\La(M)=\bigoplus_{p\geqs 0}\La_p(M)$ we denote the exterior
$\C^\infty(M)$ algebra of vector fields on $M$.\index{$\La(M)$}
For a $\g$ module $V$, we denote by $V^\g$ the set of all $\g$
invariant elements of $V$, i.e. \linebreak
$V^\g=\{v\in V\vert X.v=0 \mbox{ for all } X\in\g\}$.

\subsection{Invariant polyvector fields on a homogeneous manifold}
\label{subsect:PolyvFieldsOnHomogManif}

We show here how to construct polyvector fields on $M=G/K$ starting 
from elements of $\bw^p\g$ and $\bw^p\g/\k$.

Taking $\ol{\tht}\in\left(\bw^p\g/\k\right)^\k$, one can construct a
$p$-vector field $\la_G(\ol{\tht})$ on $G$ in the following way.
Using reductivity of $\k$, we choose a $\k$ invariant subspace $\m$
of $\g$, complement to $\k$.
Then lift $\ol{\tht}$ to $\tht\in\left(\bw^p\g\right)^\k$ and put 
$(\la_G(\ol{\tht}))_g=(L_g)_*\tht$,
where $L_g:G\to G$ is the left translation $x\mapsto gx$.
The field $\la_G(\ol{\tht})$ is left $G$ invariant.
It is also right $K$ invariant, thus it is projectable on $M=G/K$.
Set $\la_M(\tht)=\pi_*(\la_G(\tht))$\index{$\la_M$}
where $\pi$ is the natural projection $G\to G/K$.
The field $\la_M(\tht)$ on $M$ is $G$ invariant.
We shall denote the $\C^\infty(M)$ module
of all $\g$ invariant $p$-vector fields on $M$ by $\La^\g_p(M)$.
\begin{prop}\label{th:GInvVectFlds}
Any $\g$ invariant polyvector field on $M$ is of the form
$\la_M(\tht)$ for some $\tht\in\left(\bw^p\g/\k\right)^\k$.
\end{prop}
\begin{proof}
A $\g$ invariant polyvector field is fully determined by its value 
at an arbitrary point $m\in M$.
There is a natural vector space isomorphism between the quotient 
$\g/\k$ and the tangent space $\T_mM$.
This isomorphism induces an isomorphism between $\bw^p\g/\k$ and 
$\bw^p\T_mM$.
Thus any $\g$ invariant $p$-vector field is generated by some 
element of $\bw^p\g/\k$.
On the other hand, to be projectable, this element of $\bw^p\g/\k$
must be $\k$ invariant.
(See also \cite{OniszikAL:TopOfTransfGroups}, 1.4.6).
\end{proof}

In a similar way, consider an element $\psi\in\bw^p\g$ and put 
$(\rho_G(\psi))_g=(R_g)_*\psi$\index{$\rho_M$} for any $g\in G$,
where $R_g:G\to G$ is the right translation $x\mapsto xg$.
By its definition,
$\rho_G(\psi)$ is $G$ invariant from the right.
In particular, it is $K$ invariant from the right.
This makes it projectable on $M$,
denote the corresponding polyvector field by $\rho_M(\psi)$.
Note that the field $\rho_M(\psi)$ is not necessarily $\g$ invariant.

\subsection{The Schouten bracket}
\label{subsect:SchoutenBrackets}

\begin{defn}\label{def:SchoutenBrackets}
Let $\La(\a)$ be the exterior algebra of a Lie algebra $\a$.
\emph{The Schouten bracket}\index{S{\nm chouten} bracket} on 
$\La(\a)$ is a $\kk$ bilinear mapping 
$\La_p(\a)\times\La_q(\a)\to\La_{p+q-1}(\a)$ defined by
$$
\dl X_1\w\dots\w X_p,Y_1\w\dots\w Y_q \dr=
\sum_{i,j}(-1)^{i+j}[X_i,Y_j]\w X_1\w\dots\hat{X_i}\dots
\hat{Y_j}\dots\w Y_q
$$
where $[X,Y]$ is the Lie bracket on $\a$ and the notation $\hat{X_i}$ 
means that $X_i$ is omitted in the summand.
\end{defn}
In particular, putting $\a=\Vect(M)$, the {\nm Lie} algebra of vector
fields on $M$,\index{$\Vect(M)$} one obtains the bracket 
$\La_p(M)\times\La_q(M)\to\La_{p+q-1}(M)$ of polyvector fields, 
introduced by {\rf J.~A.~Schouten} in \cite{SchoutenJA:DiffKomit}.
Interpreting polyvectors on $M$ as polyderivations of the function
algebra $A=\C^\infty(M)$, one can give an equivalent definition of
their {\nm Schouten} bracket \cite{LichnerowiczA:NATO},
\cite{DeWildeLecomte:NATO}.
We give here the corresponding formulas for two particular cases
which are important for us.

Let $\xi,\eta\in\La_2(M)$, then the {\nm Schouten} bracket 
$\dl\xi,\eta\dr$ is a $3$-vector field acting on a triple
$a\ot b\ot c\in A^{\ot 3}$ as follows:
$$
\dl\xi,\eta\dr(a,b,c)=
\xi(\eta(a,b)c)+\eta(\xi(a,b),c)+\mathrm{Cycl}
$$
where Cycl denotes taking of all the cyclic permutations of 
$a,b,c$.
Using this presentation of {\nm Schouten} bracket, we prove that
\begin{equation}\label{eq:AltIsSchoutBr}
\dl\xi,\eta\dr=
3\Alt\Big(\xi\comp(\eta\ot\id)+\eta\comp(\xi\ot\id)\Big).
\end{equation}
Indeed, by using the definition of the operator $\Alt$ and the skew 
symmetry of $\xi$ and $\eta$, one obtains:
\begin{eqnarray*}
\Alt\Big(
\xi\comp(\eta\ot\id)\Big)(a,b,c)
&=&
\frac{1}{6}\Big(
\xi\left(\eta(a,b),c\right)-\xi\left(\eta(b,a),c\right)+ \\
&+&\xi\left(\eta(b,c),a\right)-\xi\left(\eta(c,b),a\right)+ \\
&+&\xi\left(\eta(c,a),b\right)-\xi\left(\eta(a,c),b\right)\Big)= \\
&=&\frac{1}{3}\Big(\xi\left(\eta(a,b),c\right)+
\xi\left(\eta(b,c),a\right)+\xi\left(\eta(c,a),b\right)\Big)= \\
&=&\frac{1}{3}\Big(\xi\left(\eta(a,b),c\right)+\mathrm{Cycl}
\Big).
\end{eqnarray*}
In the same way,
$
\Alt\Big(\eta\comp(\xi\ot\id)\Big)(a,b,c)=
\frac{1}{3}\Big(\eta\left(\xi(a,b),c\right)+\mathrm{Cycl}
\Big),
$
which proves (\ref{eq:AltIsSchoutBr}).

For $\xi\in\La_2(M)$, $\ups\in\La_3(M)$, one has
\begin{eqnarray}\label{eq:SchBrTwoWithThree}
\dl\xi,\ups\dr&=&2\Alt\Big(
\xi\comp(\ups\ot\id)-\ups\comp(\xi\ot\id\ot\id)+\\
&+&\ups\comp(\id\ot\xi\ot\id)-\ups\comp(\id\ot\id\ot\xi)
+\xi\comp(\id\ot\ups)\Big). \nonumber
\end{eqnarray}

The following proposition will be useful for cohomological
calculations.
\begin{prop}\label{th:VectFieldsOnHomManifold}
 \ \nl
\mbox{\hspace{.2in} \emph{(I)}\hspace{.2in}}
Let $\ol{\tht},\ol{\ups}\in\La^\k(\g/\k)$, 
then $\dl\la_M(\ol{\tht}),\la_M(\ol{\ups})\dr=
\la_M(\dl\ol{\tht},\ol{\ups}\dr)$.\nl
\mbox{\hspace{.17in} \emph{(II)}\hspace{.17in}}
Let $\tht,\upsilon\in\La(\g)$, 
then $\dl\rho_M(\tht),\rho_M(\upsilon)\dr=
-\rho_M(\dl\tht,\upsilon\dr)$.\nl
\mbox{\hspace{.15in} \emph{(III)}\hspace{.15in}}
Let $\tht\in\La(\g)$, $\upsilon\in\La^\k(\g)$, 
then $\dl\rho_M(\tht),\la_M(\upsilon)\dr=0$.
\end{prop}
The left-hand side of each formula is the {\nm Schouten} bracket of
polyvector fields on $M$ while the right-hand side contains the
{\nm Schouten} bracket on $\La(\g)$.

\begin{proof}
First we consider polyvector fields generated on the group $G$ and
then pass to the manifold $M=G/K$.
Using Definition~\ref{def:SchoutenBrackets}, one sees that
the {\nm Schouten} brackets in $\La(\g)$ and in $\La(G)$ defined
by the {\nm Lie} bracket in $\g$ and the commutator of
left invariant vector fields on $G$ correspondingly.
Since the {\nm Lie} bracket in $\g$ generated by the commutator of
left invariant vector fields on $G$, one has 
$\la_G([\ol{X},\ol{Y}])=[\la_G(\ol{X}),\la_G(\ol{Y})]$ for 
$\ol{X},\ol{Y}\in\g/\k$ and
$\rho_G([X,Y])=-[\rho_G(X),\rho_G(Y)]$ for $X,Y\in\g$.
Thus
$\la_G(\dl\ol{\tht},\ol{\ups}\dr)=
\dl\la_G(\ol{\tht}),\la_G(\ol{\ups})\dr$ 
for $\ol{\tht},\ol{\ups}\in\La^\k(\g/\k)$ and
$\rho_G(\dl\tht,\ups\dr)=-\dl\rho_G(\tht),\rho_G(\ups)\dr$ for
$\tht,\ups\in\La(\g)$.
This implies the second claim of the proposition immediately.
To complete the proof of the first claim, note that if $\ol{\tht}$ 
and $\ol{\ups}$ are $\k$ invariant then the same true for
$\dl\ol{\tht},\ol{\ups}\dr$.
Indeed, it is easy to check that for any $X\in\k$ one has
$X.\dl\ol{\tht},\ol{\ups}\dr=
\dl X.\ol{\tht},\ol{\ups}\dr+\dl\ol{\tht},X.\ol{\ups}\dr$.
Thus for the polyvector fields $\la_G(\ol{\tht})$ and 
$\la_G(\ol{\ups})$ projectable their {\nm Schouten} bracket
$\dl\rho_G(\tht),\rho_G(\ups)\dr$ is projectable too.
This completes the proof for (I).

The third claim is obvious since left and right invariant vector 
fields on $G$ are always commuting.
\end{proof}

\subsection{$\ph$-Poisson brackets}
\label{subsect:PhPoissBrackets}

We showed in Section~\ref{sect:InfinitesimalsForAlgDeforms}
that the infinitesimal element of an algebra deformation of 
$\C^\infty(M)$ is necessarily a bivector field on $M$.
A bivector field obeying the {\nm Jacobi} identity is a {\nm Poisson}
bracket.
We introduce here a bracket on a homogeneous manifold which is an 
analog of {\nm Poisson} bracket for the category with non-trivial 
associator $\Ph$.
Fix a non-zero element $\wtl{\ph}\in\left(\bw^3\g\right)^\g$ and put
$\ph=\la_M(\wtl{\ph})$.\index{$\ph$}
\begin{defn}\label{def:PhiPoiss}
Let $M$ be a $G$-manifold.
A skew symmetric biderivation 
$\{\cdot,\cdot\}: \C^\infty(M)^{\ot 2}\to\C^\infty(M)$ is called 
a \emph{$\ph$-Poisson bracket}\index{$\ph$-{\nm Poisson} bracket} on 
$M$ if for any $f,g,h\in\C^\infty(M)$ one has
$$
\{\{f,g\},h\}+\{\{h,f\},g\}+\{\{g,h\},f\}=\ph(f,g,h).
$$
\end{defn}
\begin{prop}\label{th:mCYBEIsJacobi}
Let $\ups\in\La_2(M)$, then
the expression $\dl\ups,\ups\dr$ is the left-hand side of the Jacobi
identity for the bracket corresponding to $\ups$.
\end{prop}
\begin{proof}
Straightforward computation.
\end{proof}

Thus $\ups\in\La_2(M)$ corresponds a {\nm $\ph$-Poisson} bracket on 
$M$ if and only if $\dl\ups,\ups\dr=\ph$.
It determines a Poisson structure if and only if
$\ph=\la_M(\wtl{\ph})=0$.
It will be shown in Section~\ref{sect:Examples} that using an 
appropriate r-matrix bivector field, one can turn any manifold with 
$\ph$-{\nm Poisson} bracket into a {\nm Poisson} manifold.
(It is crucial that the $3$-vector fields $\la_M(\wtl{\ph})$ and 
$\rho_M(\wtl{\ph})$ coincide and $G$ invariant.)
Moreover, {\rf J.~Donin, D.~Gurevich} and {\rf S.~Shnider} recently 
proved (\cite{DoninGurevicShnider:DoubleQuantiz}, Proposition~2.2.)
that the infinitesimal of any $\Ughr$ invariant quantization
is always of the form $s+r$ where $s$ is a $\ph$-{\nm Poisson} bracket,
$r=\rho_M(\wtl{r})$ and $\wtl{r}$ is a {\nm Belavin--Drinfeld}
r-matrix.

\newpage
\section{Quantization of $\ph$-Poisson manifolds}
\label{sect:PhiPoissQuantization}

In this section we give a sufficient condition for existence for a
given $\ph$-{\nm Poisson} bracket $s$ on $M$ a $G$ invariant
$\Pb$-associative deformation of the multiplication in $\C^\infty(M)$ 
with $s$ as the infinitesimal.

As above, $G$ denotes a simple {\nm Lie} algebra over $\kk$ with
the {\nm Lie} algebra, $K$ a connected closed {\nm Lie} subgroup of 
$G$, $\k$ the {\nm Lie} algebra of $K$.
By $s$ we always denote a $\ph$-{\nm Poisson} bracket on $M=G/K$ 
generated by an element $\wtl{s}\in(\bw^2\g)^\k$ as it was explained 
in Section~\ref{subsect:PolyvFieldsOnHomogManif}.

\subsection{The complex $\mathbf{(\wtl{\La}(M),\di_s)}$}
\label{subsect:ComplexLaTilde}

For  a {\nm $\ph$-Poisson} bracket $s\in\La_2^\g(M)$ set 
$\di_s(\ups)= \dl s,\ups\dr$, $\ups\in\La^\g(M)$.
\begin{prop}\label{th:SchoutCochainCmplx}
$(\La^\g(M),\di_s)$ is a cochain complex.
\end{prop}
\begin{proof}
We need to prove that $\di_s\comp\di_s=0$.
Take $\ups\in\La_p^\g(M)$, then 
$\ups=\la_M(\wtl{\ups})$ for some $\wtl{\ups}\in(\bw^p\g/\k)^\k$
(see Section~\ref{subsect:PolyvFieldsOnHomogManif}).
We denote by $\dl\ \cdot,\cdot\ \dr$ the {\nm Schouten} brackets both
on $\La^\g(M)$ and on $\La^\k(\g/\k)$.
Thus $\di_s(\ups)=\la_M(\dl\wtl{s},\wtl{\ups}\dr)$, so it suffices
to prove that $\dl\wtl{s},\dl\wtl{s},\wtl{\ups}\dr\dr=0$ for any
$\wtl{\ups}\in(\bw^p\g)^\k$.
One checks easily that 
$\dl\wtl{s},\dl\wtl{s},\wtl{\ups}\dr\dr=
\frac{1}{2}\dl\dl\wtl{s},\wtl{s}\dr,\wtl{\ups}\dr=
\frac{1}{2}\dl\wtl{\ph},\wtl{\ups}\dr$ where $\wtl{\ph}$ a $\g$ 
invariant element of $\bw^3\g$.
All this implies that $\di_s(\di_s(\ups))=\frac{1}{2}\dl\ph,\ups\dr$
where $\ph=\la_M(\wtl{\ph})$.
Note that also $\ph=\rho_M(\wtl{\ph})$ thus, by 
Proposition~\ref{th:VectFieldsOnHomManifold}(III), 
$\dl\ph,\ups\dr=0$.
\end{proof}

Now, take a $G$ invariant field $\ups=\la_M(\wtl{\ups})$
on $M$, then one has $\di_s(\ups)= \la_M(\dl\wtl{s},\wtl\ups\dr)$.
Thus the cochain complexes $(\La^\g(M),\di_s)$ and 
$(\La^\k(\g/\k),\dl\wtl{s},\ \cdot\ \dr)$ are isomorphic.

Suppose $\La(\g)$ possesses an involution $\Caut:\La(\g)\to\La(\g)$ 
such that $\Caut(\wtl{\omg})=(-1)^{p+1}\wtl{\omg}$ if 
$\wtl{\omg}\in\bw^p\g$.
Since $\Caut(\wtl{s})=-s$,the set of all elements with 
property $\Caut(\wtl{\omg})=\wtl{\omg}$ forms 
a sub-complex of $\left(\La^\k(\g/\k),\di_s\right)$ which we denote
by $\left(\La^{\k,\tht}(\g/\k),\dl\wtl{s},\cdot\dr\right)$.
The corresponding sub-complex of $(\La^\g(M),\di_s)$ we denote by 
$\left(\wtl{\La}(M),\di_s\right)$.\index{$\wtl{\La}(M)$}
\index{$\La^{\k,\tht}(\g/\k)$}

\subsection{Effect of Drinfeld associator}

We prove here that the presence of the {\nm Drinfeld} associator 
$\Pb$ of Theorem~\ref{th:QzationBelDrinfMatrix} does not affect 
substantially the cohomological construction of the deformed 
multiplication $\mh$.
The reason is basically that $\Pb$ begins with $1\ot 1\ot 1$
and has no term with $\hb$, so it does not change the infinitesimal.

Let $\mh^{n-1}$ be a $\g$ invariant $\kk\dhb$ linear mapping 
$\fAh\ot_R\fAh\to\fAh$ of the form 
$\mh^{n-1}=\mu+ \sum_{k=1}^{n-1}\hb^k\mu_k$.
(The pair $(\fAh,\mh^{n-1})$ is a not necessary associative algebra over
$\kk\dhb$.
\begin{lm}\label{th:PhDoesNotMatter}
Set $B_\Ph(\mu)=\mu\comp(\mu\ot\id)-\mu\comp(\id\ot\mu)\comp\Ph$
and suppose that $B_{\Pb}(\mh^{n-1})\equiv 0\mdhs[n]$.
Then the following congruences are valid modulo $\hb^{n+2}$:\nl
$B_{\Pb}(\mh^{n-1})\comp\Pb^{-1}\equiv
B_{\Pb}(\mh^{n-1}) $;\nl
$B_{\Pb}(\mh^{n-1})\comp(\id\ot\De\ot\id)
(\Pb)\comp(\Ph\ot\id)\equiv 
B_{\Pb}(\mh^{n-1})  $;\nl
$B_{\Pb}(\mh^{n-1})
\comp(\id\ot\mh^{n-1}\ot\id)\comp(\Pb\ot\id)\equiv
B_{\Pb}(\mh^{n-1})\comp(\id\ot\mh^{n-1}\ot\id) $;\nl
$B_{\Pb}(\mh^{n-1})\comp(\id\ot\id\ot\mh^{n-1})\comp
(\De\ot\id\ot\id)(\Pb)\equiv
B_{\Pb}(\mh^{n-1})\comp(\id\ot\id\ot\mh^{n-1}) $.
\end{lm}
\begin{proof}
The congruence $B_{\Pb}(\mh^{n-1})\equiv 0\mdhs$ implies 
\begin{equation}\label{eq:PhiAssocNOne}
B_{\Pb}(\mh^{n-1})
\equiv\hb^n\eta+\hb^{n+1}\xi
\mdh[n+2]
\end{equation}
for some $\eta,\xi\in\wtl\C^3(A;A)$.
The formal power series $\Pb$ has $1\ot 1\ot 1$ as the initial term, 
and it has no linear term.
Its formal inverse, $\Pb^{-1}$, is of the same form.
Thus
\begin{eqnarray*}
(\hb^n\eta+\hb^{n+1}\xi)\comp\Pb & \equiv 
& \hb^n\eta+\hb^{n+1}\xi  \mdh[n+2]\\
(\hb^n\eta+\hb^{n+1}\xi)\comp\Pb^{-1} & \equiv 
& \hb^n\eta+\hb^{n+1}\xi  \mdh[n+2].
\end{eqnarray*}
Hence the left hand side of (\ref{eq:PhiAssocNOne}) will also 
not be changed after taking the composition with any algebraic
expression containing $\Pb$ or $\Pb^{-1}$.
\end{proof}

Denote by $\ol\di$\  the operator $\fAh^{\ot p}\to\fAh^{\ot (p+1)}$
defined by formula (\ref{eq:DefHochschCobound}) with $\mh^{n-1}$
instead of $\mu$.
It is a polynomial in $\hb$ with the {\nm Hochschild} 
coboundary operator $\di$ (defined by $\mu$) as the initial term:
$\ol\di=
\di+\hb\di_1+\hb^2\di_2+\cdots+\hb^{n-1}\di_{n-1}$.
In what follows we need the explicit formula for $\di_1$:
\begin{equation}\label{eq:DefDOne}
\di_1\xi=s\comp(\id\ot\xi)+\sum_{k=1}^p(-1)^k
\xi\comp(\id^{\ot(k-1)}\ot s\ot\id^{\ot(p-k)}) + 
(-1)^{p+1}s\comp(\xi\ot\id).\ \ \ \ 
\end{equation}
The operator $\ol\di$\  is not a differential, since 
$\ol\di\comp\ol\di\not=0$.
However, using the fact that $\mh^{n-1}$ is associative modulo 
$\hb^2$, we prove the following lemma.
\begin{lm}\label{th:ObstructionIsHochCocycle}
Let $\mh^{n-1}$ be $\g$ invariant and let
$B_{\Pb}(\mh^{n-1})\equiv 0\mdhs$.
Then $\ol\di\  B_{\Pb}(\mh^{n-1})\equiv 0\mdhs[n+2]$. 
\end{lm}
\begin{proof}
For this proof we set the following notation:
$B=B_{\Pb}(\mh^{n-1})$, $\ol\mu=\mh^{n-1}$ and $\Pb=\Ph$.
The sign $\equiv$ will denote the congruence modulo $\hb^{n+2}$.

By the definition of $\ol\di$\ , 
\begin{eqnarray}\label{eq:CoboundOfBTilde}
\ol\di B&=&
\ol\mu\comp(\id\ot B)
-B\comp(\ol\mu\ot\id\ot\id)+\\
&+&B\comp(\id\ot\ol\mu\ot\id)
-B\comp(\id\ot\id\ot\ol\mu)
+\ol\mu\comp(B\ot\id). \nonumber
\end{eqnarray}
Compute each term of the right hand side of 
(\ref{eq:CoboundOfBTilde}).
The following two identities can be checked directly:
$
\id\ot\left(\ol\mu\comp(\ol\mu\ot\id)\right)=
(\id\ot\ol\mu)\comp(\id\ot\ol\mu\ot\id)
$
and
$
\id\ot\left(\ol\mu\comp(\id\ot\ol\mu)\comp\Ph\right)=
(\id\ot\ol\mu)\comp(\id\ot\id\ot\ol\mu)\comp(\id\ot\Ph).
$
Thus 
\begin{equation}\label{eq:PhiAssocHbNTwoFirst}
\ol\mu\comp(\id\ot B)=
\ol\mu\comp(\id\ot\ol\mu)\comp(\id\ot\ol\mu\ot\id)-
\ol\mu\comp(\id\ot\ol\mu)\comp(\id\ot\id\ot\ol\mu)\comp
(\id\ot\Ph)
\end{equation}
In analogous way, using the identities
$\left(\ol\mu\comp(\ol\mu\ot\id)\right)\ot\id=
(\ol\mu\ot\id)\comp(\ol\mu\ot\id\ot\id)$ and
$\left(\ol\mu\comp(\id\ot\ol\mu)\comp\Ph\right)\ot\id=
(\ol\mu\ot\id)\comp(\id\ot\ol\mu\ot\id)\comp(\Ph\ot\id)$,
one obtains
\begin{eqnarray}\label{eq:PhiAssocHbNTwoFifth}
\ol\mu\comp(B\ot\id)&=&
\ol\mu\comp(\ol\mu\ot\id)\comp(\ol\mu\ot\id\ot\id)-\\
&-&\ol\mu\comp(\ol\mu\ot\id)\comp(\id\ot\ol\mu\ot\id)
\comp(\Ph\ot\id) \nonumber
\end{eqnarray}

To proceed with the remaining three terms of 
(\ref{eq:CoboundOfBTilde}),
recall that $\ol\mu$ is $\g$ invariant,
i.e. $\ol x\comp\ol\mu=\ol\mu\comp\De(\ol x)$ for all
$\ol x\in\ol U$ (see Definition~\ref{def:GInvarMult}).
Since $\Ph\in\ol U^{\ot 3}$,
this implies the following three equations:
$$
\begin{array}{l}
\Ph\comp(\ol\mu\ot\id\ot\id)=
(\ol\mu\ot\id\ot\id)\comp
(\De\ot\id\ot\id)(\Ph);\\
\ \\
\Ph\comp(\id\ot\ol\mu\ot\id)=
(\id\ot\ol\mu\ot\id)\comp
(\id\ot\De\ot\id)(\Ph);\\
\ \\
\Ph\comp(\id\ot\id\ot\ol\mu)=
(\id\ot\id\ot\ol\mu)\comp
(\id\ot\id\ot\De)(\Ph).
\end{array}
$$
Using these equations, one obtains
\begin{eqnarray}\label{eq:PhiAssocHbNTwoSecond}
B\comp(\ol\mu\ot\id\ot\id)&=&
\ol\mu\comp(\ol\mu\ot\id)\comp(\ol\mu\ot\id\ot\id)- \nonumber\\
&-&\ol\mu\comp(\id\ot\ol\mu)\comp\Ph\comp
(\ol\mu\ot\id\ot\id)=\\
&=&\ol\mu\comp(\ol\mu\ot\id)\comp(\ol\mu\ot\id\ot\id)-\nonumber\\
&-&\ol\mu\comp(\id\ot\ol\mu)\comp(\ol\mu\ot\id\ot\id)
\comp(\De\ot\id\ot\id)(\Ph); \nonumber
\end{eqnarray}

\begin{eqnarray}\label{eq:PhiAssocHbNTwoThird}
B\comp(\id\ot\ol\mu\ot\id)&=&
\ol\mu\comp(\ol\mu\ot\id)\comp(\id\ot\ol\mu\ot\id)-\\
&-&\ol\mu\comp(\id\ot\ol\mu)\comp(\id\ot\ol\mu\ot\id)
\comp(\id\ot\De\ot\id)(\Ph); \nonumber
\end{eqnarray}

\begin{eqnarray}\label{eq:PhiAssocHbNTwoFourth}
B\comp(\id\ot\id\ot\ol\mu)&=&
\ol\mu\comp(\ol\mu\ot\id)\comp(\id\ot\id\ot\ol\mu)\\
&-&\ol\mu\comp(\id\ot\ol\mu)\comp(\id\ot\id\ot\ol\mu)
\comp(\id\ot\id\ot\De)(\Ph). \nonumber
\end{eqnarray}

Take the composition from the right of the both
sides of (\ref{eq:PhiAssocHbNTwoFirst})
with $(\id\ot\De\ot\id)(\Ph)\comp(\Ph\ot\id)$.
By Lemma~\ref{th:PhDoesNotMatter},
$B\comp(\id\ot\De\ot\id)(\Ph)\comp(\Ph\ot\id)\equiv B$.
Thus one obtains the left hand side of 
(\ref{eq:PhiAssocHbNTwoFirst}) unchanged,
and so
\begin{eqnarray}\label{eq:PhiAssocHbNTwoFirstPrim}
\ol\mu\comp(\id\ot B)&\equiv&
\ol\mu\comp(\id\ot\ol\mu)\comp(\id\ot\ol\mu\ot\id)\comp
(\id\ot\De\ot\id)(\Ph)\comp(\Ph\ot\id)-
\\
&-&\ol\mu\comp(\id\ot\ol\mu)\comp(\id\ot\id\ot\ol\mu)
\comp (\id\ot\Ph)\comp(\id\ot\De\ot\id)(\Ph)\comp
(\Ph\ot\id) \nonumber
\end{eqnarray}
Similarly, taking the composition from the right of the both
sides of (\ref{eq:PhiAssocHbNTwoThird})
with $\Ph\ot\id$ and the both sides of 
(\ref{eq:PhiAssocHbNTwoFourth}) with $(\De\ot\id\ot\id)(\Ph)$,
one obtains:
\begin{eqnarray}\label{eq:PhiAssocHbNTwoThirdPrim}
B\comp(\id\ot\ol\mu\ot\id)&\equiv&
\ol\mu\comp(\ol\mu\ot\id)\comp(\id\ot\ol\mu\ot\id)
\comp(\Ph\ot\id) - \\
&-&\ol\mu\comp(\id\ot\ol\mu)\comp(\id\ot\ol\mu\ot\id)
\comp(\id\ot\De\ot\id)(\Ph)
\comp(\Ph\ot\id) \nonumber
\end{eqnarray}

\begin{eqnarray}\label{eq:PhiAssocHbNTwoFourthPrim}
B\comp(\id\ot\id\ot\ol\mu)&\equiv&
\ol\mu\comp(\ol\mu\ot\id)\comp(\id\ot\id\ot\ol\mu)
\comp(\De\ot\id\ot\id)(\Ph)- \\
&-&\ol\mu\comp(\id\ot\ol\mu)\comp(\id\ot\id\ot\ol\mu)
\comp(\id\ot\id\ot\De)(\Ph)
\comp(\De\ot\id\ot\id)(\Ph) \nonumber
\end{eqnarray}

Finally, put (\ref{eq:PhiAssocHbNTwoFirstPrim}),
(\ref{eq:PhiAssocHbNTwoSecond}),
(\ref{eq:PhiAssocHbNTwoThirdPrim}),
(\ref{eq:PhiAssocHbNTwoFourthPrim}) and
(\ref{eq:PhiAssocHbNTwoFifth}) into (\ref{eq:CoboundOfBTilde}).
Using the identity
$(\ol\mu\ot\id)\comp(\id\ot\id\ot\ol\mu)=
\ol\mu\ot\ol\mu=
(\id\ot\ol\mu)\comp(\ol\mu\ot\id\ot\id)$
and the Pentagon Identity, one concludes that $\ol\di\ B\equiv 0$.
\end{proof}

\subsection{Quantization in the category with $\Pb$-associativity}

Let $\ph$ and $s$ be as above, $(\wtl\C^\bullet(M),\di\,)$ the 
local {\nm Hochschild} complex for $A=\C^\infty(M)$,
$\Pb$ as in Theorem~\ref{th:QzationBelDrinfMatrix}.
Our goal is to prove that there exists a non-commutative 
$\Pb$-associative formal deformation $(A\dl\hb\dr,\mh)$ of the 
associative commutative algebra $(A,\mu)$ with the  original 
associative multiplication $\mu$ as the initial term and $s$ as the 
infinitesimal.
Recall (Example~\ref{eg:PhAssoc}) that $\Pb$ associativity means 
$\mh\comp(\mh\ot\id)=\mh\comp(\id\ot\mh)\comp\Pb$.
This can be written as
\begin{equation}\label{eq:PhiHbAssociativity}
B_{\Pb}(\mh)=0.
\end{equation}
Recall also that $\g$ invariance means $x\comp\mh=\mh\comp\De(x)$ for 
any $x\in\Ug$ considering as a differential operator on $A\dl\hb\dr$.

Suppose that $\La(\g)$ is equipped with an involution $\Caut$.
Recall (see Section~\ref{subsect:ComplexLaTilde}) that the subset
$\wtl{\La}(M)$ of $\La(\g)$ consisting of elements whose homogeneous
components $\wtl{\omg}\in\bw^p\g$ satisfy  
$\Caut(\wtl{\omg})=(-1)^{p+1}(\wtl{\omg})$ forms a sub-complex.
\begin{defn}\label{def:PhQuatability}
We call a $\ph$-Poisson bracket $s$ 
\emph{$\Pb$-quantizable}\index{$\Pb$-quantizability} if there exists 
a $\Pb$-associative commutative $\g$ invariant multiplication 
$\mu_\hb$ on $\C^\infty(M)\dl\hb\dr$ of the form 
$\mh=\mu+\hb s+\cdots$.
\end{defn}
\begin{thm}\label{th:PhQuanizationOfSExists}
Let $s$ be a $\g$ invariant bivector field on $M$ such that 
$\dl s,s\dr=\ph$.
Suppose that $\Ha^3\left(\wtl{\La}(M),\di_s\right)=0$.
Then $s$ is $\Pb$-quantizable.
\end{thm}
This theorem is a generalization of Proposition~4 of 
\cite{DoninShnider:QSymSpaces}.
The difference of our considerations from those of 
\cite{DoninShnider:QSymSpaces} is that we fix the first order 
term for $\mh$ and we do not assume $\la_M(\ph)=0$.
Since $\Ha^3(\wtl\C^\bullet(M),\di\, )\not=0$ in general,
one can not apply the arguments of \cite{DoninShnider:QSymSpaces} 
directly to the case.
Instead, we combine these arguments with the techniques of  
{\rf O.~M.~Nero\-slav\-sky} and {\rf A.~T.~Vlassov}
\cite{NeroslavskijWlassow}.
A clear exposition of their method in the context of obstruction 
theory is given in \cite{LichnerowiczA:NATO}.

\begin{proof}
Starting with $\mh^1=\mu+\hb s$, we find subsequent 
approximations of $\mu_\hb$ by $\g$ invariant multiplications 
$$
\mu_\hb^n=\mu + s\hb + \sum_{i=2}^n\mu_i\hb^i
$$
with $\mu_i$'s obeying the Parity Convention: $\mu_{2k}\in\wtl\C_-^2(A;A)$ and 
$\mu_{2k+1}\in\wtl\C_+^2(A;A)$.

To proceed by induction on $n$, we prove first that the linear approximation,
$\mh^1=\mu+\hb s$, obeys (\ref{eq:PhiHbAssociativity}) modulo $\hb^2$, i.e. 
$$
B_{\Pb}(\mh^1)=
\mh^{1}\comp(\mh^{1}\ot\id)-
\mh^{1}\comp(\id\ot\mh^{1})\comp\Pb
\equiv 0\mdh[2]
$$
Opening the brackets and leaving the constant and linear on $\hb$ 
terms only and using the definition (\ref{eq:DefHochschCobound}) of 
the {\nm Hochschild} coboundary operator, one obtains 
$B_{\Pb}(\mh^1)=\di\mu+\hb\di s$.
Simple computation shows that $\di\mu=0$.
We prove that $\di s=0$.
Using the definition (\ref{eq:DefHochschCobound})
of operator $\di$ and the {\nm Leibniz} rule 
for each argument of $s$,
one has for any $a,b,c\in A$:
\begin{eqnarray*}
\di s(a,b,c)&=&as(b,c)-s(ab,c)+s(a,bc)-s(a,b)c=  \\
&=&as(b,c)-as(b,c)-bs(a,c)+s(a,c)b+s(a,b)c-s(a,b)c=0,
\end{eqnarray*}
where $ab=\mu(a,b)$.
Note that $\Pb$ still played no r\^ole in our arguments.

Assume we found $\mu_i$'s for all $i<n$, such that $\mu_\hb^{n-1}$
is $\Pb$-associative modulo $\hb^n$:
\begin{equation}\label{eq:PhiAssocN}
B_{\Pb}(\mh^{n-1})=
\mh^{n-1}\comp(\mh^{n-1}\ot\id)-
\mh^{n-1}\comp(\id\ot\mh^{n-1})\comp\Pb
\equiv 0 \mdh
\end{equation}
This means that
\begin{equation}\label{eq:PhiAssocNOnePrim}
B_{\Pb}(\mh^{n-1})
\equiv\hb^n\eta_n
\mdh[n+1]
\end{equation}
for some $\eta_n\in\wtl\C^3(A;A)$.
We are looking for a $\mu_n\in\wtl\C^2(A;A)$ which will cancel
$\eta_n$, 
i.e. such that $\mh^n=\mh^{n-1}+\hb^n\mu_n$ obeys 
(\ref{eq:PhiAssocN}) modulo $\hb^{n+1}$.
Opening the brackets in (\ref{eq:PhiAssocNOne}) and putting all 
the terms of degree $n$ together one obtains an explicit 
expression for the obstruction $\eta_n$:
\begin{equation}\label{eq:EtaN}
\eta_n=
\hspace{-.1in}\sum_{\begin{array}{c}
\scriptstyle
i+j=n\\
\scriptstyle 
0< i,j<n
\end{array}}
\hspace{-.1in}\mu_i\comp(\mu_j\ot\id)-
\hspace{-.1in}\sum_{\begin{array}{c}
\scriptstyle
i+j+2k=n\\
\scriptstyle 
0<i,j<n\\
\scriptstyle 
k\geqs 0
\end{array}}
\hspace{-.1in}\mu_i\comp(\id\ot\mu_j)\comp\ph_{2k},
\end{equation}
where $\ph_{2k}$ as in Theorem~\ref{th:QzationBelDrinfMatrix}.
Note that the obstruction $\eta_2$ is $\g$ invariant, for this is 
true for $\mu$ and $\mu_1=s$.
Thus we assume that $\eta_n$ is $\g$ invariant and prove that $\mu_n$
can be chosen to be $\g$ invariant as well.
\begin{lm}\label{th:EtaIsCocycle}
  \ \nl
\mbox{\hspace{.2in}\emph{(I)}\hspace{.2in}}
$\di\eta_n=0$;\nl
\mbox{\hspace{.17in}\emph{(II)}\hspace{.17in}}
$\di_s(\Alt\eta_n)=0$.
\end{lm}
\begin{proof}
Since $B_{\Pb}(\mh^{n-1})\equiv 0\mdhs$, one has
$$
B_{\Pb}(\mh^{n-1})\equiv
\hb^n\eta_n+\hb^{n+1}\xi
\mdh[n+2]
$$
for some $\xi$.
By Lemma~\ref{th:ObstructionIsHochCocycle},
$\ol\di B_{\Pb}(\mh^{n-1})\equiv 0\mdhs[n+2]$, thus
\begin{eqnarray*}
\ol\di B_{\Pb}(\mh^{n-1})&=&
(\di+\hb\di_1+\cdots)
(\hb^n\eta_n+\hb^{n+1}\xi+\cdots)=\\
&=&\hb^n\di\eta_n+
\hb^{n+1}(\di_1\eta_n+\di\xi)+\cdots
\equiv 0\mdh[n+2]
\end{eqnarray*}
where $\di$ is the {\nm Hochschild} coboundary for the algebra 
$(A,\mu)$,
and $\di_1$ is as in (\ref{eq:DefDOne}).
Equating to zero the coefficients before $\hb^n$ and $\hb^{n+1}$
respectively, one obtains two equations:
$\di\eta_n=0$ and $\di_1\eta_n+\di\xi=0$.
The former equation proves (I).
Taking the alternation of the both sides of the latter equation, one 
obtains $\Alt(\di_1\eta_n)=0$ since $\Alt(\di\xi)=0$
by Proposition~\ref{th:AltOfHochCoboundIsZero}.
Using $\di\eta_n=0$ and Theorem~\ref{th:KostHochRosenbKonc},
one concludes that $\Alt\eta_n$ is a $3$-vector field on $M$.
Using direct computation, one can prove that 
$2\Alt(\di_1\eta_n)=\dl s,\Alt \eta_n\dr$ (see formula
(\ref{eq:SchBrTwoWithThree}) and the proof for formula 
(\ref{eq:AltIsSchoutBr})).
This proves (II).
\end{proof}

In classical obstruction theory, when the third 
cohomology space is equal to zero, all obstructions vanish,
but we have $\Ha^3(\wtl\C^\bullet(M),\di\,)\not=0$.
To eliminate the obstruction $\eta_n$,
one can try to correct the previous term, $\mu_{n-1}\hb^{n-1}$.
First prove that the parity of $\eta_n$ is the opposite to the parity
of the integer $n$, i.e. 
$\eta_{2k}\in\C^3_{-}(A;A)$ and $\eta_{2k+1}\in\C^3_{+}(A;A)$.
\begin{lm}\label{th:EtaNEvenOdd}
$\eta_n(c,b,a)=(-1)^{n+1}\eta_n(a,b,c)$.
\end{lm}
\begin{proof}
By Lemma~\ref{th:PhDoesNotMatter}, 
$B_{\Pb}(\mh^{n-1})\comp\Pb^{-1}\equiv B_{\Pb}(\mh^{n-1})\mdhs[n+1]$.
This implies
\begin{equation}\label{eq:IndepOnPhi}
\begin{array}{l}
\mh^{n-1}\comp(\mh^{n-1}\ot\id)-
\mh^{n-1}\comp(\id\ot\mh^{n-1})\comp\Pb
 \equiv   \\
  \   \\
 \equiv 
\mh^{n-1}\comp(\mh^{n-1}\ot\id)\comp\Pb^{-1}-
\mh^{n-1}\comp(\id\ot\mh^{n-1})
\mdh[n+1]
\end{array}
\end{equation}

For a formal power series $\xi=\sum_{k\geqs 0}\xi_k\hb^k$,
set $\mathbf{c}_n(\xi)=\xi_n$.
Clearly, $\mathbf{c}_n(\xi)=0$ if $\xi\equiv 0\mdhs[n+1]$.
Therefore (\ref{eq:IndepOnPhi}) implies that
\begin{eqnarray}\label{eq:EtaInTwoForms}
\eta_n&=&
\mathbf{c}_n\left(\mh^{n-1}\comp(\mh^{n-1}\ot\id)-
\mh^{n-1}\comp(\id\ot\mh^{n-1})\comp\Pb\right)= \\
  &=&\mathbf{c}_n
\left(\mh^{n-1}\comp(\mh^{n-1}\ot\id)
\comp\Pb^{-1}-
\mh^{n-1}\comp(\id\ot\mh^{n-1})\right) \nonumber
\end{eqnarray}

Assume $n$ is even, then the parity of $i$ and $j$ in each term of 
(\ref{eq:EtaN}) is the same.
Consequently,
$\mu_i\left(\mu_j(a,b),c\right)=\mu_i\left(c,\mu_j(b,a)\right)$,
and thus
$$
\mh^n(\mh^n(a,b),c))=\mh^n(c,\mh^n(b,a))
$$
for any $a,b,c\in A$.
Using (\ref{eq:EtaInTwoForms}) and the property $\Pb^{321}=\Pb^{-1}$,
one obtains:
\begin{equation}\label{eq:EtaCBA}
\eta_n(c,b,a)=
\mathbf{c}_n\Big(
\mh^{n-1}\comp(\mh^{n-1}\ot\id)\comp
\Pb^{321}(c,b,a)-
\mh^{n-1}\comp(\id\ot\mh^{n-1})(c,b,a)
\Big).
\end{equation}
The second summand in (\ref{eq:EtaCBA}) is equal
to $\mh^{n-1}\comp(\mh^{n-1}\ot\id)(a,b,c)$.
To consider the first summand,
write $\Pb(a,b,c)$ as
$\Ph'(a)\ot\Ph''(b)\ot\Ph'''(c)$,
then $\Ph^{321}(c,b,a)=\Ph'''(c)\ot\Ph''(b)\ot\Ph'(a)$,
and thus
\begin{eqnarray*}
\mh^{n-1}\comp(\mh^{n-1}\ot\id)\comp
\Pb^{321}(c,b,a)&=&
\mh^{n-1}\left(
\Ph'(a),\mh^{n-1}\left(
\Ph''(b),\Ph'''(c)
\right)\right)=\\
 &=&
\mh^{n-1}\comp(\id\ot\mh^{n-1})
\comp\Pb(a,b,c).
\end{eqnarray*}
One concludes that $\eta_n(c,b,a)=-\eta_n(a,b,c)$.

Similarly, when $n$ is odd, then the indices $i$ and $j$ are of the 
opposite parity for each term in (\ref{eq:EtaN}).
Thus 
$\mh^n\left(\mh^n(a,b),c)\right)=-\mh^n\left(c,\mh^n(b,a)\right)$
which is followed by $\eta_n(c,b,a)=\eta_n(a,b,c)$.
\end{proof}

By Lemma~\ref{th:EtaIsCocycle}, $\eta_n$ is a {\nm Hochschild} 
cocycle.
Then by Theorem~\ref{th:KostHochRosenbKonc}, 
$\Alt\eta_n$ is a $3$-vector field on $M$, and
\begin{equation}\label{eq:DecomposEta}
\eta_n=\Alt\eta_n+\di\nu
\end{equation}
for some bidifferential operator $\nu$.
Note that $\Alt\eta_n$ is $\g$ invariant.

Assume $n$ is odd.
By Lemma~\ref{th:EtaNEvenOdd}, $\eta_n\in\wtl{\C}^3_+(A;A)$, i.e.
it is even, and thus $\Alt\eta_n$ is equal to zero, for this is 
an odd $3$-cochain.
Thus in this case $\eta_n=\di\nu$ for some bidifferential 
operator $\nu\in\wtl{\C}^2(A;A)$.
Proposition~\ref{th:TauIsComplexMorphism} shows that $\nu$ in 
(\ref{eq:DecomposEta}) can be chosen skew symmetric.
\begin{lm}\label{th:HochCoboundCanBeChsnInvariant}
Let $\xi=\di_H\ups$ be a $\g$ invariant Hochschild 
$p$-coboundary, $\ups\in\wtl{\C}^{p-1}(A;A)$.
Then there exists a $\g$ invariant 
$\ups^\circ\in\wtl{\C}^{p-1}(A;A)$ such that 
$\xi=\di_H\ups^\circ$.
\end{lm}
\begin{proof}
Consider the spaces $\wtl{\C}^{p}(A;A)$ as $\g$ modules.
The modules $\wtl{\C}^{p}(A;A)$ can be decomposed into direct sum of
finite dimensional highest weight spaces.
The sum of the spaces of highest weight zero is equal exactly to the
subspace $\wtl{\C}^{p}(A;A)^\g$ of $\g$ invariant cochains.
For $\wtl{\C}^p=\wtl{\C}^{p}(A;A)$ and a weight $\la$ of $\g$, 
denote by $V^\la(\wtl{\C}^p)$ the direct sum of all components of the
highest weight $\la$.
In particular, $V^0(\wtl{\C}^p)=(\wtl{\C}^p)^\g$.
It is easy to check that the {\nm Hochschild} coboundary operator 
$\di_H$ commutes with the $\g$ action.
Thus 
\begin{equation}\label{eq:DiHochRespectsHighWeights}
\di_H\left(V^\la(\wtl{\C}^{p-1})\right)\subset V^\la(\wtl{\C}^p).
\end{equation}
Decompose $\ups=\ups^{\circ}+\ups^{'}$ where 
$\ups^\circ\in V^0(\wtl{\C}^{p-1})$ and 
$\ups^{'}\in\dSum_{\la\not=0}V^\la(\wtl{\C}^{p-1})$.
Then, by (\ref{eq:DiHochRespectsHighWeights}), 
$\xi\in V^0(\wtl{\C}^{p})$ implies $\di_H\ups^{'}=0$, and thus
$\xi=\di_H\ups^{\circ}$. 
\end{proof}

\noindent
This lemma implies that $\nu$ can be chosen $\g$ invariant.
The polynomial $\mh^n=\mh^{n-1}+\hb^n\nu$ is $\Pb$-associative 
modulo $\hb^{n+1}$ and it is $\g$ invariant.
Hence for an odd $n$, $\g$ invariant term $\mu_n\hb^n$ obeying the 
parity convention can always be constructed.

Consider the case when $n$ is even.
By Lemma~\ref{th:EtaNEvenOdd}, $\eta_n$ is odd, thus $\Alt\eta_n$ is 
not necessary zero.
Proposition~\ref{th:TauIsComplexMorphism} and 
Lemma~\ref{th:HochCoboundCanBeChsnInvariant} show that $\nu$ in 
(\ref{eq:DecomposEta}) can be chosen symmetric and $\g$ invariant.
By Lemma~\ref{th:EtaIsCocycle}, $\di_s(\Alt\eta_n)=0$.
Together with $\Ha^3\left(\wtl{\La}(M),\di_s\right)=0$, this implies 
that $\Alt\eta_n=\di_s\ze$ for some $\ze\in\wtl{\La}_2(M)$.

Put 
$$
\tmh^{n-1}=\mh^{n-1}-\frac{3}{2}\ze\hb^{n-1}.
$$
Since $\ze$ is skew symmetric, $\tmh^{n-1}$ obeys the parity 
convention.
Thus for the obstruction corresponding to $\tmh^{n-1}$, i.e. a 
bidifferential operator $\wtl\eta_n$ such that
$B_{\Pb}(\tmh^{n-1})\equiv\hb^n\wtl\eta_n\mdhs[n+1]$, the statements 
of Lemma~\ref{th:EtaIsCocycle} (II) and Lemma~\ref{th:EtaNEvenOdd} 
are still valid, that is, $\wtl\eta_n(c,b,a)=-\wtl\eta_n(a,b,c)$
and $\di\wtl\eta_n=0$.
Direct computation shows that
\begin{equation}\label{eq:TildEtaN}
\wtl\eta_n=
\eta_n+\frac{3}{2}\hb^{n-1}\di \ze-
\frac{3}{2}\Big(
\ze\comp(s\ot\id)+s\comp(\ze\ot\id)-
\ze\comp(\id\ot s)-s\comp(\id\ot\ze)
\Big)
\end{equation}
The {\nm Hochschild} coboundary of any biderivation is equal to 
zero, hence $\di\ze=0$. 
Using the skew symmetry of $s$ and $\ze$, one checks directly that
$$
\Alt\Big(\ze\comp(\id\ot s)+s\comp(\id\ot\ze)\Big)=
-\Alt\Big(\ze\comp(s\ot\id)+s\comp(\ze\ot\id)\Big).
$$
Thus taking the alternation of the both sides of (\ref{eq:TildEtaN}) 
and using formula (\ref{eq:AltIsSchoutBr}),
Section~\ref{subsect:SchoutenBrackets}, one obtains 
$\Alt\wtl\eta_n=\Alt\eta_n-\dl s,\ze\dr=\tht-\di_s\ze=\tht-\tht=0$.
By Theorem~\ref{th:KostHochRosenbKonc}, this implies that 
$\wtl\eta_n=\di\ups$ for some bidifferential operator $\ups$.
As it was mentioned above, $\ups$ can be chosen to be symmetric and
$\g$ invariant.
Finally, put $\mu_n=\ups$ and $\mh^n=\wtl{\mh}^{n-1}+\hb^{n}\ups$.
It is left to prove that $B_{\Pb}(\mh^n)\equiv 0 \mdhs[n+1]$.
Using the definition of $\wtl{\eta}_n$ and equality
$\wtl{\eta}_n=\di\ups$, one has
$B_{\Pb}(\wtl{\mh}^{n-1})\equiv\hb^n\di\ups \mdhs[n+1]$.
Therefore
\begin{eqnarray*}
B_{\Pb}(\mh^n)&=&
B_{\Pb}(\wtl{\mh}^{n-1})+\hb^n\Big(
\mu\comp(\ups\ot\id)+\ups\comp(\mu\ot\id)\Big)-\\
&\ &-\hb^n\Big(
\mu\comp(\id\ot\ups)+\ups\comp(\id\ot\mu)\Big)\comp
\Pb\equiv \\
&\equiv&B_{\Pb}(\wtl{\mh}^{n-1})+\hb^n\Big(
\mu\comp(\ups\ot\id)+\ups\comp(\mu\ot\id)-\\
&\ &-\mu\comp(\id\ot\ups)-\ups\comp(\id\ot\mu)\Big)\equiv\\
&\equiv&\di\ups-\di\ups\equiv 0 \mdh[n+1],
\end{eqnarray*}
where the congruence $\hb^n\mu\comp(\id\ot\ups)\comp
\Pb\equiv\hb^n\mu\comp(\id\ot\ups)\mdhs[n+1]$ was used 
(see Lemma~\ref{th:PhDoesNotMatter}).
Theorem~\ref{th:PhQuanizationOfSExists} is proven.
\end{proof}

\newpage
\section{A class of homogeneous manifolds \\ 
with quantizable Poisson brackets}
\label{sect:Examples}

In this section we introduce a class of homogeneous manifolds, $\Mal$,
which closely related to manifolds appearing in the problem of 
classification of quotients of $G$ by a reductive subgroup of 
maximal rank.
We prove that all the manifolds $\Mal$ posses $G$ invariant 
$\ph$-{\nm Poisson} brackets and present their explicit forms.
It turns out that these brackets are essentially unique.
Any $\ph$-{\nm Poisson} bracket $s$ determines a cochain complex.
By computing the corresponding cohomologies of this complex, we prove 
that $s$ can be quantized in such a way
that after quantization we obtain a $G$ invariant $\Pb$-associative 
multiplication.
Due to results of Section~\ref{sect:MonCategDrAlgs}, this implies that
any {\nm Poisson} bracket on $\Mal$ of the form $s+r$ can be
quantized invariantly with respect to the quantum group $\Ughr$ 
action.

\subsection{Poisson brackets generated by Belavin--Drinfeld 
$r$-matrices}

Let $G$ be a connected simple {\nm Lie} group over $\kk$, $\g$ the
{\nm Lie} algebra of $G$, $K$ a connected {\nm Lie} subgroup of $G$,
$\k$ the {\nm Lie} algebra of $K$.
Denote by  $M$ the homogeneous $G$ manifold $G/K$.
Recall that the vector space $\m=\g/\k$ is isomorphic to the tangent
space to $M$ at the point fixed by $K$.
In Section~\ref{subsect:PhPoissBrackets} we introduced the notion of
$\ph$-{\nm Poisson} bracket.
Here we explain how to pass from a $\ph$-{\nm Poisson} bracket to 
a usual {\nm Poisson} bracket.

Let $s$ be a $\ph$-{\nm Poisson} bracket on $M$.
By Proposition~\ref{th:GInvVectFlds}, $s=\la_M(\wtl{s})$
for some $\wtl{s}\in\left(\bw^2\m\right)^\k$ satisfying 
$\dl\wtl{s},\wtl{s}\dr=\wtl{\ph}$, where $\wtl{\ph}$ is a non-zero
$\g$ invariant element of $\bw^3\g$.
Take a {\nm Belavin---Drinfeld} r-matrix $\wtl{r}$ with 
$\dl\wtl{r},\wtl{r}\dr=\wtl{\ph}$ for the same $\wtl{\ph}$, and 
put $r=\rho_M(\wtl{r})$.
Then, by Proposition~\ref{th:VectFieldsOnHomManifold}, one has
$\dl s,s\dr=\ph$, $\dl r,r\dr=-\ph$ and $\dl s,r\dr=0$, thus
$\dl s+r,s+r\dr=0$.

As in Section~\ref{subsect:ExistenceOfFHbPhiHb}, we denote by 
$\Ughr$ the {\nm Etingof--Kazhdan} quantum group determined by $r$.
Recall that $\Pb$-quantizability of $s$ means that there exists a 
$\Pb$-associative $\g$ invariant multiplication $\mu_\hb$ on 
$\C^\infty(M)\dl\hb\dr$ of the form $\wtl{\mh}=\mu+\hb s+\cdots$ 
where $\Pb$ is given by Theorem~\ref{th:QzationBelDrinfMatrix}.
\begin{thm}\label{th:QzationOfRPlusS}
Let $M$ be a smooth manifold with the above $\ph$-bracket $s$ and 
{\nm Belavin---Drinfeld} bivector field.
Then \nl
\mbox{\hspace{.2in} \emph{(I)}\hspace{.2in}}
$s+r$ is a Poisson bracket on $M$;\nl
\mbox{\hspace{.17in} \emph{(II)}\hspace{.17in}}
If $s$ is $\Pb$-quantizable, then there exists an associative 
multiplication $\mh$ on $\C^\infty(M)\dhb$ of the form 
$\mh=\mu+\hb(s+r)+\cdots$ which is invariant under the action of the 
quantum group $\Ughr$.
\end{thm}
\begin{proof}
(I).\ 
It was mentioned above that $\dl s+r,s+r\dr=0$.
By Proposition~\ref{th:mCYBEIsJacobi}, this means that the sum $s+r$
obeys the {\nm Jacobi} identity.

(II).\ 
Suppose that there exists a $\Pb$-associative multiplication of the 
form $\wtl{\mu}_\hb=\mu+\hb s+\cdots$.
By Theorem~\ref{th:QzationBelDrinfMatrix}, there exists an invertible
power series $\Fb=\id^{\ot 2}+\hb r+\cdots\in\U^{\ot 2}\g\dhb$ such 
that the composition $\mh=\wtl{\mu}_\hb\comp\Fb$ is a (strictly) 
associative multiplication.
Expanding the composition in powers of $\hb$, one obtains 
$\mh=\mu+\hb(s+r)+\cdots$.
(See also Remark~\ref{nb:MultiplsInEquivCategs}.)
\end{proof}

\subsection{Homogeneous manifolds related to regular subalgebras of 
$\g$}
\label{subsect:RegularSubAlgs}

Fix a {\nm Cartan} subalgebra $\Car$ of the simple {\nm Lie} algebra 
$\g$, denote by $\Wurz$ the corresponding root system and fix a set 
of simple roots for $\Wurz$.

In what follows, all tensor products and dimensions are taken over 
$\kk$.
\begin{defn}
We call a Lie subalgebra $\k$ \emph{regular}\index{L{\nm ie}
subalgebra!regular}\index{regular Lie subalgebra} if it is reductive 
and contains a Cartan subalgebra.
\end{defn}
Since all {\nm Cartan} subalgebras of $\g$ are conjugate, we can
consider only those regular subalgebras which contain the fixed
{\nm Cartan} subalgebra $\Car$.
Choose a subset $\UWurz\subset\Wurz$ and denote by $\Ga(\Wurz)$ and 
$\Ga(\UWurz)$ the $\z$ lattices generated by $\Wurz$ and $\UWurz$ 
respectively.
Set $\Wurz_\UWurz=\Ga(\UWurz)\cap\Wurz$ and denote by $\k$ the 
{\nm Lie} subalgebra in $\g$ of the form
$$
\k=\Car\dsum\left(\dSum_{\be\in\Wurz_\UWurz}\g^\be\right)
$$
where $\g^\be$\index{$\g^\be$} is the one dimensional  root space in 
$\g$ corresponding to $\be\in\Wurz$, and $\dsum$ denotes the direct 
sum of vector spaces.
Then $\k$ is regular, and any regular {\nm Lie} subalgebra of $\g$ 
appears in this way, see \cite{GorbazewicOnishchikVinverg:LieAlgIII},
Chapter~6, \S~1.
Denote by $K$ the {\nm Lie} subgroup of $G$ corresponding to the
{\nm Lie} algebra $\k$.
Since $K$ corresponds to a regular {\nm Lie} algebra, it is 
connected and closed.
Set $M=G/K$ and\index{$\m$}
$$
\m=\dSum_{\be\in\Wurz\setminus\Wurz_{\UWurz}}\g^\be.
$$
It is a $\k$ module, and it is easy to see that $\m$ is the 
orthogonal complement to $\k$ with respect to the {\nm Killing} form 
$(\cdot,\cdot)$.
Thus one can identify $\m$ with the quotient $\g/\k$ and treat it as
the tangent space to the homogeneous $M$ at the point fixed by $K$.

Denote by $\ol{\Wurz}$\index{$\ol{\Wurz}$} the image without zero of
$\Wurz$ under the canonical epimorphism 
$\Ga(\Wurz)\to\Ga(\Wurz)/\Ga(\UWurz)$.
We denote the image of $\be\in\Wurz$ in $\ol{\Wurz}$ by $\clbe$,
\index{$\clbe$} and call the elements of $\ol{\Wurz}$
\emph{quasi-roots}\index{quasi-root}.
Quasi-roots are convenient labels for some important irreducible 
representations of {\nm Lie} algebra $\k$.
Namely, put
$$
\m_{\clbe}=\dSum_{\ga\in\clbe}\g^\ga.
$$\index{$\m_{\clbe}$}
To prove that all $\k$ modules $\m_{\clbe}$ are simple, we will use
the following lemma which is proven in 
\cite{DoninGurevicShnider:DoubleQuantiz}.
\begin{lm}\label{th:EquivRootsAreConnected}
Let $\clbe=\ol{\be^{'}}$, then there exist
$\al_1,\dots,\al_k\in\Wurz_\UWurz$ such that
$\be+\al_1+\al_2+\cdots+\al_k=\be^{'}$ and 
$\be+\al_1+\al_2+\cdots+\al_i$ is a root for every $i\leqs k$.
\end{lm}
\begin{proof}
We use the following fact about root systems of simple {\nm Lie}
algebras (see \cite{SerreJP:AlgLieSS}, Chap. IV, Prop.~3).
If $(\al,\be)>0$ for some roots $\al$, $\be$ then $\al-\be$ is a
root as well.
The equality $\ol{\be^{'}}=\clbe$ means that 
${\be^{'}}=\be+\ga_1+\cdots+\ga_m$ for some (not necessary
different) roots $\ga_i\in\UWurz$.
If $(\be^{'},\be)>0$ then $\be^{'}-\be$ is a root in 
${\Wurz}_{\UWurz}$, and everything is done.
In the opposite case, $(\be^{'},\be)\leqs 0$, one uses induction by
$m$.
Namely, since $(\be^{'},\be^{'})>0$, there exists a root $\ga_i$ such
that $(\be^{'},\ga_i)>0$.
Changing labels, one can assume $(\be^{'},\ga_m)>0$.
Thus $\be^{'}-\ga_m=\be+\ga_1+\cdots\ga_{m-1}$ is a root, and
$\ol{\be^{'}-\ga_m}=\clbe$.
By the induction hypothesis, there exist 
$\al_1,\dots,\al_{m-1}\in\Wurz_\UWurz$ such that
$\be+\al_1+\al_2+\cdots+\al_{m-1}=\be^{'}-\ga_m$ and 
$\be+\al_1+\al_2+\cdots+\al_i$ is a root for every $i\leqs m-1$.
Taking $\al_m=\ga_m$, one completes the proof.
\end{proof}

\begin{cor}\label{th:mBarBeIsIrred}
For any $\clbe\in\ol{\Wurz}$, $\m_{\clbe}$ is an irreducible $\k$ 
module.
\end{cor}
\begin{proof}
Chose a base $\{X_\al\}$ of weight vectors for $\m_{\clbe}$ such 
that $X_\al\in\g^\al$, and take $\be$, $\be^{'}$ and $\al_i$ as in 
Lemma~\ref{th:EquivRootsAreConnected}.
Then $X_{\al_i}\in\k$, and
$X_{\be^{'}}=\ad X_{\al_k}\cdots\ad X_{\al_2}\ad X_{\al_1}(X_\be)$.
This implies that any element of a basis of weight vectors of 
$\m_{\clbe}$ can be mapped into another arbitrary element with the 
help of a composition of operators from $\ad\k$.
Thus all the weight spaces are in same irreducible component.
\end{proof}

\noindent
The following lemma is also proven in 
\cite{DoninGurevicShnider:DoubleQuantiz}.
\begin{lm}\label{th:SumHasReptatives}
Let $\ol{\be},\ol{\be}_1,\dots,\ol{\be}_m\in\ol{\Wurz}$ such that 
$\ol{\be}=\sum_{i=1}^m\ol{\be}_i$.
Then there exist roots $\ga,\ga_i\in\Wurz$ such that
$\ga\in\ol{\be}$, $\ga_i\in\ol{\be}_i$ and $\sum_{i=1}^m\ga_i=\ga$.
\end{lm}
\begin{proof}
The equality $\clbe=\sum_{i=1}^m\ol{\be}_i$ means that 
$\be=\sum_{i=1}^m {\be}_i+\sum_{j=1}^n\al_j$ for some 
$\al_1,\al_2,\dots,\al_n\in\Wurz_\UWurz$.
Let $(\be,\al_j)>0$ for some $\al_j$, then $\be^{'}=\be-\al_j$ is a 
root representing $\clbe$.
This procedure can be repeated if necessary several times.
Thus one can assume that $(\be,\al_j)\leqs 0$ for all $j$'s.
We use the induction on $m$ to find the representatives $\ga_i$.
For $m=1$, there is nothing to prove.
Assume that the lemma is proven for all sums of the form
$\ol{\be}=\sum_{i=1}^{m-1}\ol{\be}_i$.
Let $\be=\sum_{i=1}^m {\be}_i+\sum_{j=1}^n\al_j$ and 
$(\be,\al_j)\leqs 0$ for all $j$'s.
Since $(\be,\be)=\sum_{i=1}^m(\be,\be_i)+\sum_{j=1}^n(\be,\al_j)>0$,
there exists $i$, say $i=m$, such that $(\be,\be_m)>0$.
Therefore $\be-\be_m$ is a root.
Put $\be^{''}=\be-\be_m$, then 
$\ol{\be^{''}}=\sum_{i=1}^{m-1}\ol{\be}_i$, i.e. $\ol{\be^{''}}$ is 
a sum of $m-1$ quasi-roots, and the induction hypothesis applies.
\end{proof}

\begin{cor}\label{th:mBarBemBarGaIsmBarBePlusBarGa}
$[\m_{\ol{\be}_1},\m_{\ol{\be}_2}]=\m_{\ol{\be}_1+\ol{\be}_2}$
\end{cor}
\begin{proof}
The inclusion $[\m_{\ol{\be}_1},\m_{\ol{\be}_2}]\subset
\m_{\ol{\be}_1+\ol{\be}_2}$ is obvious.
We prove that $[\m_{\ol{\be}_1},\m_{\ol{\be}_2}]\supset
\m_{\ol{\be}_1+\ol{\be}_2}$.
By Lemma~\ref{th:SumHasReptatives}, $\clbe[1]+\clbe[2]\in\ol{\Wurz}$
implies that there exist representatives $\ga_1\in\clbe[1]$ and 
$\ga_2\in\clbe[2]$ such that $\ga_1+\ga_2\in\clbe[1]+\clbe[2]$ is a
root.
Then the space $[\g^{\ga_1},\g^{\ga_2}]$ is non-zero, and it is
contained in $\m_{\clbe[1]+\clbe[2]}$.
By Corollary~\ref{th:mBarBeIsIrred}, the latter is irreducible.
This proves the claim.
\end{proof}

\noindent
Since $\k$ is reductive, the $\k$ module $\m$ can be decomposed 
into direct sum of irreducible modules.
\begin{lm}\label{th:mBarAlIrredInTensProd}
Let $\ol{\be}_1+\ol{\be}_2\in\ol{\Wurz}$, then 
$\m_{\ol{\be}_1+\ol{\be}_2}$ is a multiplicity free irreducible  
component of $\k$ module $\m_{\ol{\be}_1}\ot\m_{\ol{\be}_2}$.
\end{lm}
\begin{proof}
First we prove that the $\k$ module $\m_{\ol{\be}_1+\ol{\be}_2}$ 
appears as a component of $\m_{\ol{\be}_1}\ot\m_{\ol{\be}_2}$.
Consider the mapping $\cL:\m_{\clbe[1]}\ot\m_{\clbe[2]}\to
[\m_{\ol{\be}_1},\m_{\ol{\be}_2}]$,  $X\ot Y\mapsto[X,Y]$.
It is a $\k$ module homomorphism.
By Corollary~\ref{th:mBarBemBarGaIsmBarBePlusBarGa}, $\cL$ has image
$\m_{{\clbe}_1+{\clbe}_2}$.
The {\nm Lie} algebra $\k$ is reductive, thus the $\k$ module
$\m_{\clbe[1]}\ot\m_{\clbe[2]}$ is decomposed into direct sum of 
irreducible components.
Among them, there exists an irreducible sub-module, say
$\mathfrak{p}$, such that the restriction $\cL\vert\mathfrak{p}$
is non-zero.
By {\nm Schur}'s Lemma, it is an isomorphism of $\k$ modules.

The fact that $\m_{\ol{\be}_1+\ol{\be}_2}$ is multiplicity free
follows from \cite{Bourbaki:GroupLieAlgLie78}, Chap.~{VIII}, \S 9,
Ex. 14).
\end{proof}

Recall that if $V$ is a module over an arbitrary {\nm Lie} algebra
$\a$, the dual space $V^*=\Hom_\kk(V,\kk)$ is endowed with an $\a$
module structure: $(X\xi)(v)=-\xi(Xv)$ for $v\in V$, $\xi\in V^*$.
\begin{lm}\label{th:MinusBetaDualToBeta}
For any $\be\in\Wurz$,
$\k$ modules $\m_{-\ol{\be}}$ and $\m^{*}_{\ol{\be}}$ are isomorphic.
\end{lm}
\begin{proof}
One has $(\g^\be,\g^\ga)=0$ if and only if $\be+\ga\not=0$.
Hence the restriction of the {\nm Killing} form $(\cdot,\cdot)$ to 
$\m_{\ol{\be}}\ot\m_{-\ol{\be}}$ is non-degenerate, since
$\m_{\ol{\be}}\ot\m_{-\ol{\be}}=\sum\g^{\ga_1}\ot\g^{\ga_2}$.
The restriction is $\k$ invariant, thus it defines an isomorphism 
$\m_{-\ol{\be}}\cong\m^{*}_{\ol{\be}}$.
\end{proof}
\begin{cor}\label{th:TwoComponentIsOfDimOne}
For any $\clbe\in\ol{\Wurz}$, one has 
$\dim\left(\m_{\ol{\be}}\ot\m_{-\ol{\be}}\right)^\k=1$.
\end{cor}
\begin{proof}
First note that if $U$ and $V$ are modules over $\k$, then the $\k$
modules $(V^*\ot U)^\k$ and $\Hom_{\k}(V,U)$ are isomorphic.
Indeed take $\xi\ot u\in(V^*\ot U)^\k$ and set $L(v)=\xi(v)u$ for 
any $v\in V$.
The operator $L$ is $\k$ linear since $X(\xi\ot u)=0$ for any 
$X\in\k$.

By Lemma~\ref{th:MinusBetaDualToBeta}, $\m_{-\ol{\be}}$ and 
$\m^*_{\ol{\be}}$ are isomorphic, thus the $\k$ module
$\left(\m_{\ol{\be}}\ot\m_{-\ol{\be}}\right)^{\k}$ is isomorphic to  
$\End_{\k}(\m_{\ol{\be}})$.
By Corollary~\ref{th:mBarBeIsIrred}, $\m_{\ol{\be}}$ is
irreducible, thus one can apply {\nm Schur}'s Lemma and conclude
that all elements of $\End_\k(\m_{\ol{\be}})$ are multiples of the 
unity operator.
\end{proof}

\begin{cor}\label{th:BeMinusBeZero}
Let $\clbe$ be a quasi-root such that $\clbe=-\clbe$, then 
$\left(\bw^2\m_{\clbe}\right)^\k=0$.
\end{cor}
\begin{proof}
By Lemma~\ref{th:MinusBetaDualToBeta},
$\left(\m_{\clbe}\ot\m_{\clbe}\right)^\k$ is
isomorphic as $\k$ module to $\End_\k(\m_{\clbe})$.
Thus $\dim\left(\m_{\clbe}\ot\m_{\clbe}\right)^\k=1$.
However, $\left(\m_{\clbe}\ot\m_{\clbe}\right)^\k$ contains a 
non-trivial element, that is, the {\nm Killing} form which is 
symmetric.
Hence there is no skew symmetric elements in 
$\left(\m_{\clbe}\ot\m_{\clbe}\right)^\k$.
\end{proof}

In what follows we need the existence of one well-known base for
the {\nm Lie} algebra $\g$.
\begin{prop}\label{th:PropertsOfNAlBe}
For a semi-simple Lie algebra $\g$ there exists a base 
$\{X_\al\}_{\al\in\Wurz}$, $\{H_\be\}_{\be\in\Pi}$ with 
$(X_\al,X_{-\al})=1$ such that corresponding structural constants
$N_{\al,\be}$ have the following properties:\nl
\mbox{\hspace{.2in} \emph{(I)}\hspace{.2in}}
 $N_{\al\be}\not=0$  if and only if $\al+\be\in\Wurz$;\nl
\mbox{\hspace{.17in} \emph{(II)}\hspace{.17in}}
$N_{\al\be}=N_{\be\ga}=N_{\ga\al}$ for $\al+\be+\ga=0$;\nl
\mbox{\hspace{.14in} \emph{(III)}\hspace{.14in}}
 $N_{\be\al}=-N_{\al\be}$;\nl
\mbox{\hspace{.15in} \emph{(IV)}\hspace{.15in}}
$N_{-\al,-\be}=-N_{\al\be}$.
\end{prop}\index{$N_{\al\be}$}
\begin{proof}
See \cite{HelgasonS:DiffGeomLieAlgSymSpaces}, Lemma~III.5.1 and
Theorem~III.5.5.
\end{proof}

We are interested in $G$ invariant polyvector fields on $M$.
They are of the form 
$\la_M(\wtl{\ups})$, $\wtl{\ups}\in\left(\bw^p\m\right)^{\k}$ (see
Section~\ref{subsect:PolyvFieldsOnHomogManif}).
We describe the spaces $\left(\bw^p\m\right)^{\k}$ for $p=1,2,3$.
Obviously, there are only trivial $\k$ invariant vector fields on 
$M$ since $\k$ contains a {\nm Cartan} subalgebra.
The following lemma describes invariant $2$-vector fields on $M$. 
\begin{lm}\label{th:InvarBivectors}
Any element of $\left(\bw^2\m\right)^\k$ is of the form
$$
\wtl{\ups}=\sum_{\clal\in\ol{\Wurz}}c_{\clal}
\left(\sum_{\be\in\clal}X_\be\w X_{-\be}\right)
$$
where $c_{-\clal}=-c_{\clal}$.
\end{lm}
\begin{proof}
We prove first that any element of this form is $\k$ invariant.
It suffices to check that for any $\ga\in\Wurz_\UWurz$ and any 
quasi-root $\clal\in\ol{\Wurz}$ one has $X_\ga.\wtl{\eta}=0$ where
$$
\wtl{\eta}=\sum_{\be\in\clal}X_\be\w X_{-\be}.
$$
Choose $\be\in\clal$.
If $\be+\ga$ and $\ga-\be$ are not roots,
then $X_\ga.(X_\be\w X_{-\be})$ does not contribute to $X_\ga.\wtl{\eta}$.
If $\be+\ga$ is a root, then $X_\ga.\wtl{\eta}$ contains the terms
$N_{\ga\be}X_{\be+\ga}\w X_{-\be}$ and 
$N_{\ga,-(\be+\ga)}X_{\be+\ga}\w X_{-\be}$ which cancel one another,
see Proposition~\ref{th:PropertsOfNAlBe}.
In the same way one considers the case when $\ga-\be$ is a root.

Note that any expression of the form $\sum_{\be\in\clal}c_\be X_\be\w X_{-\be}$
is invariant with respect to the {\nm Cartan} subalgebra $\Car$.
We prove that invariance by $\k$ implies $c_\be=c_{\be^{'}}$ if
$\clbe=\ol{\be^{'}}$.
Set
$$
\wtl{\psi}=\sum_{\be\in\clal}c_\be X_\be\w X_{-\be}
$$
By Lemma~\ref{th:EquivRootsAreConnected}, if $\be\in\clal$ then 
$\be=\al+\ze_1+\cdots\ze_k$ for some
$\ze_i\in\Wurz_\UWurz$, and $\al+\ze_1+\cdots\ze_j$ is a root for any $j\leqs k$.
Let $\be=\al+\ze_1+\cdots+\ze_{i-1}$, $\ze=\ze_i$, $\be^{'}=\be+\ze$ and consider 
the element $X_\ze.\wtl{\psi}$.
Since $\be+\ze$ is a root, $X_\ze.\wtl{\psi}$ contains the term
$X_{\be+\ze}\w X_{-\be}$.
Direct computations show that its coefficient is equal to
$c_\be N_{\ze\be}+c_{\be+\ze} N_{\ze,-(\be+\ze)}$.
Using the properties of $N_{\ze\be}$, one reduces this coefficient
to the form $(c_\be -c_{\be+\ze})N_{\ze\be}$.
Thus $X_\ze.\wtl{\psi}=0$ implies $c_{\be+\ze}=c_\be$.
\end{proof}

Consider the space $\left(\bw^3\m\right)^{\k}$ of $\k$ invariant
$3$-vector fields on $M$.
First, note that both $\left(\bw^3\m\right)^{\k}$ and any its
subspace of the form
$\left(\m_{\clal}\ot\m_{\clbe}\ot\m_{\ol{\ga}}\right)^{\k}$
are invariant under the action of {\nm Cartan} subalgebra $\Car$.
This implies that any of their elements must be of weight zero, i.e. 
it is a linear combination of monomials $X_\al\ot X_\be\ot X_{\ga}$ 
with $\al+\be+\ga=0$.
In particular, if $\ol{\ga}\not=-\clal-\clbe$ then 
$\left(\m_{\clal}\ot\m_{\clbe}\ot\m_{\ol{\ga}}\right)^{\k}
=0$.
Indeed, it is easy to see that if $\ol{\ga}\not=-\clal-\clbe$ then 
$\al+\be+\ga\not=0$ for any $\al\in\clal$, $\be\in\clbe$ and
$\ga\in\ol{\ga}$.
\begin{lm}\label{th:TripleTensProdIsOfDimOne}
The dimension of 
$\left(\m_{\clal}\ot\m_{\clbe}\ot\m_{-\clal-\clbe}\right)^{\k}$ is 
one.
\end{lm}
\begin{proof}
Set $V=\m_{\clal}\ot\m_{\clbe}$ and consider this space with the 
usual $\k$ module structure.
By Lemma~\ref{th:MinusBetaDualToBeta}, $\m_{-(\clal+\clbe)}$ is
isomorphic as $\k$ module to $\m^*_{\clal+\clbe}$.
Thus the $\k$ modules $\left(V\ot\m_{-(\clal+\clbe)}\right)^\k$ and 
$\Hom_{\k}(\m_{\clal+\clbe},V)$ are isomorphic (see proof for 
Corollary~\ref{th:TwoComponentIsOfDimOne}).
Since $\k$ is reductive, $V$ is a direct sum of simple $\k$ modules.
By Lemma~\ref{th:mBarAlIrredInTensProd}, the module 
$\m_{\clal+\clbe}$ is simple and multiplicity free in 
$\m_{\clal}\ot\m_{\clbe}$.
Thus $\Hom_{\k}(\m_{\clal+\clbe},V)$ and $\End_{\k}(\m_{\clal+\clbe})$ 
are isomorphic, and $\dim\End_{\k}(\m_{\clal+\clbe})=1$.
\end{proof}

\begin{lm}\label{th:DimLaThreeGeneral}
The dimension of $\left(\bw^3\m\right)^\k$ is equal to the number of 
unordered pairs $\{\clal,\clbe\}$ such that 
$\clal+\clbe\in\ol{\Wurz}$.
\end{lm}
\begin{proof}
\noindent
The dimension of $\left(\bw^3\m\right)^\k$ is 
equal to the number of distinct subspaces of the form 
$\left(\m_{\clal}\ot\m_{\clbe}\ot\m_{-\clal-\clbe}\right)^{\k}$.
\end{proof}

\noindent
Now we describe the elements of $\left(\bw^2\m\right)^\k$
which are $\ph$-{\nm Poisson} brackets.
\begin{lm}\label{th:EqnForC}
$\dl \wtl{\ups},\wtl{\ups}\dr=\ka^2\wtl{\ph}$ if and only if the 
coefficients $c_{\clbe}$ from Lemma~\ref{th:InvarBivectors} satisfy 
the following condition: if $\clal+\clbe\in\ol{\Wurz}$ then
$$
c_{\clal+\clbe}=\frac{c_{\clal}c_{\clbe}+\ka^2}{c_{\clal}+c_{\clbe}}.
$$
\end{lm}
\begin{proof}
Let $\wtl{\omg}=\sum_{\al}b_\al X_\al\w X_{-\al}$ and
$\wtl{\ups}=\sum_{\al}c_\al X_\al\w X_{-\al}$ be elements of
$\left(\bw^2\m\right)^\k$.
We compute the coefficient before monomial
$X_{\al+\be}\w X_{-\al}\w X_{-\be}$ in $\dl\wtl{\omg},\wtl{\ups}\dr$.
This monomial appears six times, with the coefficients 
$N_{\al\be}c_{\al}b_{\be}$, $N_{-(\al+\be),\be}c_{\al+\be}b_{\be}$, 
$-N_{\be,-(\al+\be)}c_{\be}b_{\al+\be}$, 
$N_{\al,-(\al+\be)}c_{\al}b_{\al+\be}$, 
$-N_{\be\al}c_{\be}b_{\al}$ and \linebreak
$-N_{-(\al+\be),\al}c_{\al+\be}b_{\al}$. 
Using the properties of $N_{\al,\be}$ 
(Proposition~\ref{th:PropertsOfNAlBe}), one obtains the sum of the
above coefficients:
\begin{equation}\label{eq:SchBrackOfBivects}
N_{\al\be}\cdot(c_\al b_\be-c_{\al+\be} b_\be-c_\be b_{\al+\be}
-c_\al b_{\al+\be}+c_\be b_\al-c_{\al+\be} b_\al).
\end{equation}
The element $\wtl{\ph}$ is of the form 
\begin{equation}\label{eq:InvarEltOfMThree}
\wtl{\ph}=
\sum_{\begin{array}{c}
\scriptstyle
\al\in{\clal}\\
\scriptstyle 
\be\in{\clbe}
\end{array}}\hspace{-.1in}
N_{\al\be}X_\al\w X_\be\w X_{-(\al+\be)}.
\end{equation}
Thus the coefficient before $X_{\al+\be}\w X_{-\al}\w X_{-\be}$ in 
$\ka^2\wtl{\ph}$ is $-\ka^2 N_{\al\be}$.
 and, by replacing in 
$(\ref{eq:SchBrackOfBivects})$ $b_\al$'s with $c_\al$'s, one 
completes the proof.
(See also \cite{KhoroshkinRadulRubtsov:FamPoissStructHermSymSpac}).
\end{proof}

\subsection{The manifold $M_{l\al}$ and its quantization}
\label{subsect:ManifoldMal}

We preserve the notation of Section~\ref{subsect:RegularSubAlgs}.
Fix a simple root $\al$ and a positive integer $l$ not greater than
the multiplicity of $\al$ in the highest root in $\Wurz$.
Consider the set $\Wal$ of all roots in $\Wurz$ whose coefficients
before $\al$ is divisible by $l$ 
(in Section~\ref{subsect:RegularSubAlgs} we used the notation 
$\Wurz_\UWurz$).
We denote by $\ol{\Wurz}$ the set of quasi-roots for $\mal$, i.e. 
the image of $\Wal$ in $\Ga(\Wurz)/\Ga(\Wal)$, 
by $\Mal$\index{$\Mal$} the 
homogeneous manifold whose stabilizer $K$ is generated by 
$\Ga(\Wal)\cap\Wurz$. 
Clearly, $\ol{\Wurz}=\{\clal,2\clal,\ldots,(l-1)\clal\}$.
It follows from the classification of semi-simple {\nm Lie} 
algebras over $\kk$, that $l\leqs 6$. 

The significance of the manifolds $\Mal$ is in fact that their
quantization can be considered as the first step to resolving the
following more general problem.
Let $G$ be a simple connected {\nm Lie} group over $\kk$, $\g$ its
{\nm Lie} algebra, $\k\subset \g$ a reductive {\nm Lie} subalgebra, 
$K$ the corresponding {\nm Lie} subgroup of $G$.
{\nm E.~B.~Dynkin} \cite{DynkinE:SSSubAlgsInSSLieAlgs} has proven 
that the homogeneous manifold $M=G/K$ can be obtained by taking 
consequent quotients of direct products of manifolds $\Mal$.

Recall that the quotient $\m=\g/\k$ is isomorphic to the tangent 
space to $\Mal$ at the point fixed by $K$.
We calculate the dimensions of the cohomology spaces which 
figured in Theorem~\ref{th:PhQuanizationOfSExists}.

For a given $\varkappa\in\rk$, denote by $[\varkappa]$ the largest 
integer $n$ such that $n\leqs\varkappa$.
\begin{lm}\label{th:DimLaTwo}
Let $\m=\g/\k$ be the tangent space of $\Mal$, then
$\dim\left(\bw^2\m\right)^{\k}=\left[\frac{l-1}{2}\right]$. 
\end{lm}
\begin{proof}
According to Corollary~\ref{th:BeMinusBeZero} and 
Lemma~\ref{th:InvarBivectors}, this dimension is equal to the number
of quasi-roots $\clal$ such that $\clal\not=-\clal$.
The $\k$ invariance reduces that to the number of unordered pairs
$\{\clal,-\clal\}$ which is equal to $\left[\frac{l-1}{2}\right]$.
\end{proof}

Note that all elements of $\left(\bw^2\m\right)^{\k}$ are skew 
invariant with respect to the {\nm Cartan} involution defined by
\begin{equation}\label{eq:DefCartAut}
\Caut:X_\al\mapsto -X_{-\al}
\end{equation}
The image of the coboundary
map $\dl\wtl{s},\cdot\dr: 
\left(\bw^2\m\right)^{\k}\to\left(\bw^3\m\right)^{\k}$ consists of
$\tht$ invariant elements,
because, if $\wtl{\ups}\in\left(\bw^p\m\right)^{\k}$ obeys 
$\Caut(\wtl{\ups})=(-1)^{p+1}\wtl{\ups}$ then the element
$\dl\wtl{s},\wtl{\ups}\dr$ is of the opposite parity:
$\Caut(\dl\wtl{s},\wtl{\ups}\dr)=(-1)^p\dl\wtl{s},\wtl{\ups}\dr$.
Therefore one can consider the sub-complex of $\Caut$ invariant 
$p$-vector fields for $p$ odd and skew $\Caut$ invariant $p$-vector
field for $p$ even.
Recall that we denote this sub-complex by 
$\left(\bw^p\m\right)^{\k,\Caut}$.
Note that the element $\wtl{\ph}$ given by formula
(\ref{eq:InvarEltOfMThree})  belongs to 
$\left(\bw^3\m\right)^{\k,\Caut}$.

Now, we compute the dimension of $\left(\bw^3\m\right)^{\k,\Caut}$.
For quasi-roots $i\clal$, $j\clal$ such that $j\clal\not=-i\clal$, 
denote by $\vf(i,j)$ the subspace of 
$\left(\bw^3\m\right)^{\k,\Caut}$ generated by the image of $\kk$ 
linear mapping
\begin{eqnarray}\label{eq:DefOfVAlBeBar}
\left(\mal[i]\ot\mal[j]\ot\mal[{-(i+j)}]\right)^{\k}
&\to&\left(\bw^3\m\right)^{\k,\Caut}\\
X_\be\ot X_\ga\ot X_{-\be-\ga}&\mapsto& 
X_\be\w X_\ga\w X_{-\be-\ga}-
X_{-\be}\w X_{-\ga}\w X_{\be+\ga}.\nonumber
\end{eqnarray}
It is easy to see that
$$
\vf(i,j)=\vf(j,i)=\vf(i,l-i-j)=\vf(j,l-i-j),
$$
thus only spaces $\vf(i,j)$ with $i\leqs j<\frac{l}{2}$ should be 
taken into consideration.
Hence one has
\begin{equation}\label{eq:DecompOfLaThreeM}
\left(\bw^3\m\right)^{\k,\Caut}=
\dSum_{1\leqs i\leqs j<\frac{l}{2}}\vf(i,j).
\end{equation}
\begin{lm}\label{th:DimVAlBeIsOne}
The subspace $\vf(i,j)$ has dimension one.
\end{lm}
\begin{proof}
By Lemma~\ref{th:TripleTensProdIsOfDimOne},
$\dim(\mal[i]\ot\mal[j]\ot\mal[{-(i+j)}])^\k=1$.
On the other hand, the mapping (\ref{eq:DefOfVAlBeBar}) is
non-zero, since there exist $\be\in i\clal$ and $\ga\in j\clal$ such
that $\be+\ga$ is a root.
Thus the image $\vf(i,j)$ of that mapping is of dimension one.
\end{proof}

\begin{lm}\label{th:DimLaThree}
The dimension of $\left(\bw^3\m\right)^{\k,\Caut}$ is equal to the number
of its subspaces $\vf(i,j)$ with $i\leqs j<\frac{l}{2}$.
\end{lm}
\begin{proof}
Direct consequence of the formula (\ref{eq:DecompOfLaThreeM}) and
Lemma~\ref{th:DimVAlBeIsOne}.
See, also, Lemma~\ref{th:DimLaThreeGeneral}.
\end{proof}

We are interested with the third cohomology space of the 
sub-complex $\wtl{\La}(\Mal)$ of $\g$ and $\Caut$ invariant 
polyvector fields on $M$, see
Theorem~\ref{th:PhQuanizationOfSExists}.
The dimension of $\wtl{\La}_3(\Mal)$ is equal to 
$\dim\left(\bw^3\m\right)^{\k,\Caut}$.
By Lemma~\ref{th:DimLaThree}, the latter dimension is equal to the
number of pairs $(i,j)$ with $1\leqs i\leqs j<\frac{l}{2}$.
For $l=2$ there are no pairs $(i,j)$ satisfying $1\leqs i\leqs j<1$,
thus $\dim\wtl{\La}_3(\Mal[2])=0$.
For $l=3$ and $l=4$ one has one subspace in 
$\left(\bw^3\m\right)^{\k,\Caut}$, it is $\vf(1,1)$, thus
$\dim\wtl{\La}_3(\Mal[3])=\dim\wtl{\La}_3(\Mal[4])=1$.
For $l=5$ there are two subspaces in 
$\left(\bw^3\m\right)^{\k,\Caut}$, $\vf(1,1)$ and $\vf(1,2)$, thus
$\dim\wtl{\La}_3(\Mal[5])=2$.
For $l=6$ there are three subspaces in
$\left(\bw^3\m\right)^{\k,\Caut}$, $\vf(1,1)$, $\vf(1,2)$ and 
$\vf(2,2)$, thus $\dim\wtl{\La}_3(\Mal[6])=3$.

All calculated dimensions are presented in the table on 
page~\pageref{tab:Mal}.
Note that the dimension of space $\left(\bw^3\m\right)^{\k,\Caut}$ 
is equal to the number of $3$-partitions of the integer $l$, and 
that the dimension of space $\left(\bw^2\m\right)^{\k,\Caut}$ 
is equal to the number of $2$-partitions of $l$ with non-equal
components.
\begin{thm}\label{th:PhPoissBrOnMal}
 \ \nl
\mbox{\hspace{.15in}\emph{(I)}\hspace{.15in}}
Any manifold $M=\Mal$, $2\leqs l\leqs 6$, possesses a $\ph$-Poisson 
bracket.
\nl
\mbox{\hspace{.13in}\emph{(II)}\hspace{.1in}}
For any $\ph$-Poisson bracket $s$ on $M=\Mal$, $l\geqs 2$, the
cohomology spaces $\Ha^2\left(\wtl{\La}(M),\di_s\,\right)$ and
$\Ha^3\left(\wtl{\La}(M),\di_s\,\right)$ are trivial.
\end{thm}
\begin{proof}
We consider each case $l=2,\dots,6$ separately.
First prove that for $M=\Mal[2]$ one has 
$\left(\bw^3\m\right)^\Car=0$.
Indeed, take an element 
$X_{\be_1}\w X_{\be_2}\w X_{\be_3}\in\bw^3 \m$, 
$\be_1,\be_2,\be_3\not\in\Wal[2]$.
If $\be_1 + \be_2 + \be_3 \in\Ga(\Wal[2])$ then at least one of the
$\be_i$'s contains the root $\al$ with an even coefficient, and
therefore the element $X_{\be_1}\w X_{\be_2}\w X_{\be_3}$ is 
actually equal to zero.

In particular, the $3$-vector field $\ph$ is equal to zero on $M$. 
The second consequence from $\left(\bw^3\m\right)^\Car=0$ is that
$\Ha^3\left(\wtl{\La}(\Mal[2]),\di_s\,\right)=0$.

Let $M=\Mal[3]$ or $\Mal[4]$, then, by Lemma~\ref{th:DimLaTwo}
(see also the table on page~\pageref{tab:Mal}), 
$\dim\wtl{\La}_2(M)=1$.
Thus $\left(\bw^2\m\right)^{\k}$ contains a non-trivial element 
$\wtl{s}=c_{\clal}\sum_{\be\in\clal}X_\be\bw X_{-\be}$.
Direct calculation shows that $\dl\wtl{s},\wtl{s}\dr$ is 
proportional to the $\g$ invariant element $\wtl{\ph}$,
(\ref{eq:InvarEltOfMThree}), thus $s=\la(\wtl{s})$ is a 
$\ph$-{\nm Poisson} bracket on $M$.
The vector space $\left(\bw^3\m\right)^{\k,\Caut}$ is also of
dimension one (see table on p.~\pageref{tab:Mal}).
It is generated by the element $\wtl{\ph}$.
Since $\di_s(s)=\ph$, one has 
$\Ha^3\left(\wtl{\La}(M),\di_s\,\right)=0$.
On the hand, $\di_s:\wtl{\La}_2(M)\to\wtl{\La}_3(M)$ is a
non-trivial linear operator from one dimensional vector space to
another one dimensional vector space.
Thus it is non-degenerate, which proves that
$\Ha^2\left(\wtl{\La}(M),\di_s\,\right)=0$.

Let $M=\Mal[5]$, then, by Lemma~\ref{th:DimLaTwo}, 
$\dim\left(\bw^2\m\right)^{\k,\Caut}=2$.
Take a non-trivial element
$$
\wtl{s}=c_\clal\sum_{\be\in\clal}X_\be\w X_{-\be}+
c_{2\clal}\sum_{\be\in 2\clal}X_\be\w X_{-\be}\in
\left(\bw^2\m\right)^\k
$$
and prove that the multiples $c_\clal$ and $c_{2\clal}$ can be 
chosen in such a way that $\dl\wtl{s},\wtl{s}\dr=\wtl{\ph}$.
It follows from Lemma~\ref{th:EqnForC} and the identities 
$\clal+\clal=2\clal$ and $2\clal+2\clal=-\clal$ that such the 
coefficients $c_{\clal}$ and $c_{2\clal}$ should satisfy the system 
of equations
\begin{equation}
\left\{
\begin{array}{l}
5(c_{\clal})^4+10(c_{\clal})^2\ka^2+\ka^4=0 \\ 
\   \\  
c_{2\clal}=\frac{(c_{\clal})^2+\ka^2}{2c_{\clal}}. 
\end{array}
\right. \label{eq:CoefsCMFive}
\end{equation}
This system has two solutions (see table on p.~\pageref{tab:Mal}).
\begin{lm}\label{th:DTwoOnMFiveInjective}
Let $M=\Mal[5]$.
Then the coboundary operator 
$\ds:\left(\bw^2\m\right)^{\k}\to
\left(\bw^3\m\right)^{\k,\Caut}$ is injective.
\end{lm}
\begin{proof}
Any element $\wtl{\ups}$ of $\left(\bw^2\m\right)^{\k}$ is of 
the form
$$
\wtl{\ups}=b_\clal\sum_{\be\in\clal}X_\be\w X_{-\be}+
b_{2\clal}\sum_{\be\in 2\clal}X_\be\w X_{-\be}.
$$
Let $\ds \left(\wtl{\ups}\right)=0$, then the identities 
$\clal+\clal=2\clal$ and  $2\clal+2\clal=-\clal$ and
formula (\ref{eq:SchBrackOfBivects}) imply the following system of 
equations for $b_\clal$ and $b_{2\clal}$:
$$
\left\{
\begin{array}{l}
b_\clal(c_\clal+c_{2\clal})=c_\clal b_{2\clal}\\
\ \\
b_{2\clal}(c_\clal-c_{2\clal})=c_{2\clal} b_{\clal},
\end{array}
\right.
$$
where $c_\clal$ and $c_{2\clal}$ are solutions for 
(\ref{eq:CoefsCMFive}).
From the condition for the existence of a non-trivial to these
homogeneous equations and formula $(\ref{eq:CoefsCMFive})$
one comes to the following inconsistent system of equations:
$$
\left\{
\begin{array}{l}
5c_\clal^4+10c_{\clal}^2\ka^2+\ka^4=0\\
\ \\
-c_\clal^4+4c_{\clal}^2\ka^2+\ka^4=0
\end{array}
\right.
$$
Therefore the only possibility is $b_{\clal}=b_{2\clal}=0$.
\end{proof}

\noindent
The immediate consequence of this lemma is that
$\Ha^2\left(\wtl{\La}(\Mal[5]),\di_s\,\right)=0$.
On the other hand, together with the equality
$\dim\left(\bw^2\m\right)^{\k}=\dim\left(\bw^3\m\right)^{\k,\Caut}$,
it implies that  $\Ha^3\left(\wtl{\La}(\Mal[5]),\di_s\,\right)=
\Ha^3\left(\bigoplus_{p\geqs 0}(\bw^p\m)^{\k,\Caut},\ds\,\right)=0$.

Let $M=\Mal[6]$ then the table on page~\pageref{tab:Mal} shows that 
$\dim\left(\bw^2\m\right)^{\k,\Caut}=2$.
There exists a $\g$ invariant $\ph$-Poisson bracket on $M$ 
generated by an element $\wtl{s}$ of the form
$$
\wtl{s}=c_\clal\sum_{\be\in\clal}X_\be\w X_{-\be}+
c_{2\clal}\sum_{\be\in 2\clal}X_\be\w X_{-\be}
$$
with coefficients satisfying a certain equation which we obtain now.
It follows from Lemma~\ref{th:EqnForC} and the identities 
$2\clal+2\clal=-2\clal$, $\clal+\clal=2\clal$,
$\clal+2\clal=3\clal$ and $c_{3\clal}=0$ that the coefficients 
$c_{\clal}$ and $c_{2\clal}$ should satisfy the equation
\begin{equation}\label{eq:CoefsCMSix}
c_\clal^2-9 c_{2\clal}^2=0.
\end{equation}
This equation has one non-trivial solution
(see table on p.~\pageref{tab:Mal}).
\begin{lm}\label{th:DTwoOnMSixInjective}
Let $M=\Mal[6]$.
Then the coboundary operator $\ds:\left(\bw^2\m\right)^{\k}\to
\left(\bw^3\m\right)^{\k,\Caut}$ is injective.
\end{lm}
\begin{proof}
Any element $\wtl{\ups}$ of $\left(\bw^2\m\right)^{\k}$ is of the 
form
$$
\wtl{\ups}=b_\clal\sum_{\be\in\clal}X_\be\w X_{-\be}+
b_{2\clal}\sum_{\be\in 2\clal}X_\be\w X_{-\be}.
$$
Let $\ds \left(\wtl{\ups}\right)=0$, then the identities 
$2\clal+2\clal=-2\clal$ and  $\clal+\clal=2\clal$ and equations
(\ref{eq:SchBrackOfBivects}) and $(\ref{eq:CoefsCMSix})$ imply 
$b_{\clal}=b_{2\clal}=0$.
\end{proof}

\noindent
Thus $\Ha^2\left(\wtl{\La}(\Mal[6]),\di_s\,\right)=0$.
\begin{lm}\label{th:DThreeOnMSixIsNotZero}
Let $M=\Mal[6]$, then there exists
$\wtl{\ups}\in\left(\bw^3\m\right)^{\k,\Caut}$ such that
$\ds \left(\wtl{\ups}\right)\not=0$.
\end{lm}
\begin{proof}
Note that there is only one root system with the highest root
containing coefficient $6$ in the simple root decomposition.
This is the root system $E_8$, and the corresponding simple root
is the following:

\begin{picture}(300,70)(-120,-20)

\put(-40,8){$\al=\Big($}
\multiput(0,20)(20,0){7}{\circle{3}}
\multiput(1.5,20)(20,0){6}{\line(1,0){17}}
\put(40,0){\circle{3}}
\put(40,18.5){\line(0,-1){17}}
\put(130,8){$\Big)$}

\put(-3,25){$0$}
\put(17,25){$0$}
\put(37,25){$1$}
\put(57,25){$0$}
\put(77,25){$0$}
\put(97,25){$0$}
\put(117,25){$0$}

\put(37,-15){$0$}

\end{picture}

\noindent
Thus we can set $\g$ as a {\nm Lie} algebra of type $E_8$ and the 
simple root $\al$ as determined by the above weighted {\nm Dynkin} graph.

Take the element
$$
\wtl{\ups}=\sum_{\be,\ga\in\clal}N_{\be\ga}X_\be\w X_\ga\w
X_{-(\be+\ga)}
$$
and prove that $\ds\left(\wtl{\ups}\right)\not=0$.
Note that $\wtl{\ups}$ is the projection of $\wtl{\ph}$ 
onto subspace $\vf(1,1)\subset\left(\bw^3\m\right)^{\k,\Caut}$.

Choose roots $\be,\ \ga,\ \ep,\ \ze\in E_8$ satisfying the 
properties:
$\be,\ga\in 2\clal$, $\ep\in 3\clal$, $\ze\in 5\clal$,
$\be+\ga+\ep+\ze=0$; $\be+\ga$, $\be+\ze$ and $\ga+\ep$ are roots; 
$\ga+\ze$ and $\be+\ep$ are not roots.
One can take for example the following set of roots:

\begin{picture}(300,70)(-20,0)

\put(-40,8){$\be=\Big($}
\multiput(0,20)(20,0){7}{\circle{3}}
\multiput(1.5,20)(20,0){6}{\line(1,0){17}}
\put(40,0){\circle{3}}
\put(40,18.5){\line(0,-1){17}}
\put(130,8){$\Big)$}

\put(-3,25){$0$}
\put(17,25){$1$}
\put(37,25){$2$}
\put(57,25){$1$}
\put(77,25){$1$}
\put(97,25){$0$}
\put(117,25){$0$}

\put(37,-15){$1$}

\put(180,8){$\ga=\Big($}
\multiput(220,20)(20,0){7}{\circle{3}}
\multiput(221.5,20)(20,0){6}{\line(1,0){17}}
\put(260,0){\circle{3}}
\put(260,18.5){\line(0,-1){17}}
\put(350,8){$\Big)$}

\put(217,25){$1$}
\put(237,25){$2$}
\put(257,25){$2$}
\put(277,25){$2$}
\put(297,25){$2$}
\put(317,25){$2$}
\put(337,25){$1$}

\put(257,-15){$1$}

\end{picture}

\vspace{.3in}
\begin{picture}(300,70)(-24,-15)

\put(-47,8){$\ep=-\Big($}
\multiput(0,20)(20,0){7}{\circle{3}}
\multiput(1.5,20)(20,0){6}{\line(1,0){17}}
\put(40,0){\circle{3}}
\put(40,18.5){\line(0,-1){17}}
\put(130,8){$\Big)$}

\put(-3,25){$1$}
\put(17,25){$2$}
\put(37,25){$3$}
\put(57,25){$3$}
\put(77,25){$3$}
\put(97,25){$2$}
\put(117,25){$1$}

\put(37,-15){$2$}

\put(173,8){$\ze=-\Big($}
\multiput(220,20)(20,0){7}{\circle{3}}
\multiput(221.5,20)(20,0){6}{\line(1,0){17}}
\put(260,0){\circle{3}}
\put(260,18.5){\line(0,-1){17}}
\put(350,8){$\Big)$}

\put(217,25){$0$}
\put(237,25){$1$}
\put(257,25){$1$}
\put(277,25){$0$}
\put(297,25){$0$}
\put(317,25){$0$}
\put(337,25){$0$}

\put(257,-15){$0$}

\end{picture}

\vs
The monomial $X_\be\w X_\ga\w X_\ep\w X_\ze$ is contained 
in the image $\ds\,(\wtl{\ups})$.
We prove that it is present in $\ds\,(\wtl{\ups})$ with a non-zero
coefficient.

We introduce the following temporary definition.
We say that a root vector $X_\be$ is \emph{ of type $i$} if the root
$\be$ contains the simple root $\al$ with multiplicity $i$.
Obviously, $-6<i<6$.
It is clear how to extend this definition to any exterior monomial.
For instance, $X_\be\w X_\ga\w X_\ep\w X_\ze$ for the above roots is 
of type $(2,2,-3,-1)$.

The vector space $\vf(1,1)$ consists of certain sums of exterior
$3$-monomials.
The type of each monomial is a triple of integers ranging from $-5$ to
$5$.
Since $\vf(1,1)$ is invariant by the {\nm Cartan} subalgebra $\Car$, the 
sum of integers in each triple is equal to zero.
The triples satisfying these two conditions are $(1,1,-2)$, $(-1,-1,2)$, 
$(1,-5,4)$ and $(-1,5,-4)$.

We describe all possible $3$-monomials of elements in 
$\vf(1,1)$ and all possible $2$-monomials $X_\si\w X_{-\si}$ such that
their {\nm Schouten} bracket contains the given $4$-monomial 
$X_\be\w X_\ga\w X_\ep\w X_\ze$.
First note that all monomials in the {\nm Schouten} bracket of monomials 
of types $(i,j,k)$ and $(d,-d)$ are of types $(*,*,*,d)$ or 
$(*,*,*,-d)$.
Since we want to obtain a monomial of type $(2,2,-3,-1)$, this means that
$d=1,2$ or $3$.
The case $d=3$ is excluded by Corollary~\ref{th:BeMinusBeZero}.
Consider the case $d=1$.
We want to obtain a monomial of type $(2,2,-3,-1)$ by taking the 
{\nm Schouten} bracket of monomials of types $(a,b,c)$, $a+b+c=0$, and 
$(1,-1)$.
Since $[X_\al,X_\be]=N_{\al\be}X_{\al+\be}$, one concludes that each of
$a,\ b$ or $c$ is equal to either $2$ or $-3$.
Thus either $(a,b,c)=(1,2,-3)$ or $(a,b,c)=(2,2,-4)$.
However, as it was pointed out above, there is no element of $\vf(1,1)$  
containing monomials of these two types.
Let $d=2$, then the same arguments show that $a$, $b$ and $c$ are taken
from the set of integers $2$, $-3$, $-1$.
From the above list of possible types for $\vf(1,1)$, it can be seen that
only the triple $(a,b,c)=(2,-1,-1)$ fits.

Thus the monomial $X_\be\w X_\ga\w X_\ep\w X_\ze$ can be obtained as
the {\nm Schouten} bracket of monomials of types $(2,-1,-1)$ and 
$(2,-2)$.
This implies that the latter monomial is equal either to 
$c_{\clal}X_\be\w X_{-\be}$ or to $c_{2\clal}X_\ga\w X_{-\ga}$.
Since $\wtl{s}\not=0$, one has $c_{\clal}\not=0$ (see table on 
p.~\pageref{tab:Mal}).
The former monomial is either of the form 
$-6N_{\ga\ze}X_\ga\w X_{\be+\ep}\w X_\ze$ or 
$-6N_{\be\ze}X_\be\w X_{\ga+\ep}\w X_\ze$ correspondingly.
However, $\ga+\ze$ is not a root, thus the only way to obtain the
monomial $X_\be\w X_\ga\w X_\ep\w X_\ze$ is to take the 
{\nm Schouten} bracket 
$\dl -6N_{\be\ze}X_\be\w X_{\ga+\ep}\w X_\ze,
c_{2\clal}X_\ga\w X_{-\ga}\dr$.
The corresponding coefficient is $-6N_{\be\ze}c_{2\clal}$, and it 
differs from zero since $\be+\ze$ is a root.
\end{proof}

Now, we prove that 
$\Ha^3\left(\wtl{\La}(\Mal[6]),\di_s\,\right)=0$.
Indeed, the complex $\left(\wtl{\La}(\Mal[6]),\di_s\,\right)$ is 
isomorphic to the complex 
$\left(\bigoplus_{p\geqs 0}(\bw^p\m)^{\k,\Caut},\wtl{\di}\right)$.
Lemma~\ref{th:DTwoOnMSixInjective} implies that 
$$
\dim\B^3\left(\wtl{\La}(\Mal[6]),\di_s\,\right)=
\dim\wtl{\La}_2(\Mal[6])=2.
$$
On the other hand, it follows from 
Lemma~\ref{th:DThreeOnMSixIsNotZero} that 
$$
2\leqs\dim\Z^3\left(\wtl{\La}(\Mal[6]),\di_s\,\right)=
\dim \left(\bw^3\m\right)^{\k,\Caut}-\dim\left(\im\wtl{\di}^3\right)
\leqs 2.
$$
Therefore 
$$
\dim\B^3\left(\wtl{\La}(\Mal[6]),\di_s\,\right)=
\dim\Z^3\left(\wtl{\La}(\Mal[6]),\di_s\,\right).
$$
Theorem~\ref{th:PhPoissBrOnMal} is proven. 
\end{proof}

Note that $\Ha^2\left(\wtl{\La}(\Mal),\di_s\,\right)$ and 
$\Ha^3\left(\wtl{\La}(\Mal),\di_s\,\right)$ coincide with the
corresponding topological cohomologies (see \cite{BorelA:CYBE}).

In the following table we put together the computed dimensions and
give the explicit formulas for coefficients of $\ph$-{\nm Poisson}
brackets:
\vspace{.25in}

\begin{tabular}{|c||l|c|c|}
\hline
\ & \ & \ &\ \\
\ \ \ \ $l$\ \ \ \  &
\hspace{.1in} $\ph$-{\nm Poisson} brackets\hspace{.1in} & 
\ $\dim\left(\bw^2\m\right)^{\k}$ \ &
\ $\dim\left(\bw^3\m\right)^{\k,\Caut}$ \ \\
\ & \ & \ &\ \\
\hline\hline
\ & \ & \ &\ \\
2 & \hspace{.75in}  0 & 0 & 0 \\ 
\ & \ & \ & \ \\ \hline 
\ & \ & \ &\ \\
3 & $c_{\clal}=\pm\frac{i}{\sqrt{3}}\ka$ & 1 & 1 \\ 
\ & \ & \ &\ \\ \hline
\ & \ & \ &\ \\
4 & $c_{\clal}=\pm i\ka$ & 1 & 1 \\
\ & \ & \ &\ \\ \hline
\ & \ & \ &\ \\
\ & $c_{\clal}=\frac{i}{\sqrt[4]{5}}
\left(\sqrt{5}\pm 2\right)^\frac{1}{2}\ka$
& \ & \ \\
5 & $c_{2\clal}=\pm\frac{i}{\sqrt[4]{5}}
\left(\sqrt{5}\pm 2\right)^{-\frac{1}{2}}\ka$ &
2 & 2 \\
\ & \mbox{(the signs are consistent)} &\ &\ \\
\ & \ & \ &\ \\ \hline
\ & \ & \ &\ \\
\ & $c_{\clal}=\pm i\sqrt{3}\ka$ & \ & \ \\
6 & $c_{2\clal}=\pm\frac{i}{\sqrt{3}}\ka$
& 2 & 3  \\
\ & \mbox{(the signs are consistent)} &\ &\ \\
\ & \ & \ &\ \\ \hline
\end{tabular}\label{tab:Mal}

\vspace{.25in}
\noindent
Note that for $l=2,3,4,6$ there exists a $\ph$-{\nm Poisson} bracket on 
$\Mal$ unique up to a scalar multiple, and for $l=5$ there are two
such brackets.

\vs
Combining Theorems~\ref{th:PhPoissBrOnMal}, 
\ref{th:QzationOfRPlusS} and \ref{th:PhQuanizationOfSExists},
one obtains the main result of the present work:
\begin{thm}
Let $r$ be a bivector field on $\Mal$ generated by a 
Belavin--Drinfeld classical r-matrix.
Then there exists a $\g$ invariant $\ph$-Poisson bracket $s$ on
$\Mal$ such that $s+r$ is a Poisson bracket, and this bracket has a 
$\U_\hb(\g,{r})$ invariant quantization.
\end{thm}
Note that despite the fact that the bracket $s+r$ is not $\g$ 
invariant, its quantization is invariant under the quantum group 
$\Ughr$ action.

\newpage
\bibliography{BibMath}

\newpage
\ \ 
\newpage
\printindex

\end{document}